\documentclass[12pt]{article}
\usepackage{amsmath,amsfonts,amssymb,amsthm,amscd}
\usepackage{a4wide}
\title{\bf A Lindemann-Weierstrass theorem for semiabelian
varieties over function fields
}
\author{Daniel Bertrand \\Institut de Math. de Jussieu
\and Anand
Pillay   \thanks {Supported by a Marie Curie Chair and an EPSRC grant}
 \\University of Leeds}
\date{October 19, 2008}
\newtheorem{Theorem}{Theorem}[section]
\newtheorem{Proposition}[Theorem]{Proposition}
\newtheorem{Definition}[Theorem]{Definition}
\newtheorem{Remark}[Theorem]{Remark}
\newtheorem{Lemma}[Theorem]{Lemma}
\newtheorem{Corollary}[Theorem]{Corollary}
\newtheorem{Fact}[Theorem]{Fact}

\newcommand{\C}{\mathbb C}

\newcommand{\Q}{\mathbb Q}
\newcommand{\Z}{\mathbb Z}

\newcommand{\G}{\mathbb G}
\newcommand{\pf}{\noindent{\em Proof. }}

\def\mapright#1{\smash{\mathop{\longrightarrow}\limits^{#1}}}
\begin{document}
\maketitle

\begin{abstract} We prove an analogue of the Lindemann-Weierstrass theorem (that the exponentials
of a $\Q$-linearly independent set of algebraic numbers are algebraically independent), replacing $\Q^{alg}$ by
$\C(t)^{alg}$, and $\G_{m}^{n}$ by a
semiabelian variety over $\C(t)^{alg}$. Both the formulations of our results and the methods are 
differential algebraic in nature. 
\end{abstract}

\tableofcontents

\vfil \eject

\section{Introduction.}

\subsection{Statement of the results}

Let $G$ be a commutative algebraic group defined over the algebraic closure $\Q^{alg}$ of $\Q$ in $\C$, let $LG$ be its Lie algebra, and  let  $exp_G: LG(\C) \rightarrow G(\C)$ be the exponential map on the Lie group $G^{an}$ deduced from $G$ after extension of scalars to $\C$. Let further $x$ be a point in $LG(\Q^{alg})$, and assume  that

\medskip
\noindent
$({\bf HX})_{\Q^{alg}}$: for any proper algebraic subgroup $H/\Q^{alg}$ of $G$, $x \notin LH(\Q^{alg})$.

\medskip
The classical   Lindemann-Weierstrass theorem then states that if $G$ is an algebraic torus, the field of definition   of  the point  $y = exp_G(x) \in G(\C)$ satisfies:
$$ tr.deg.(\Q^{alg}(y)/\Q^{alg}) = dim(G).$$
One may wonder under which conditions on $x$  such a result extends to more general groups $G$. However, the answer is known (and with the same hypothesis)  only when $G$ is  isogenous to a power of an elliptic curve with complex multiplication (Philippon, W\"ustholz, cf.
\cite{Encyclopedia}, Theorem 6.25).  

\medskip

In the present paper, we study the analogous problem where $\Q^{alg}$ is replaced by the algebraic closure $K^{alg}$  of the function field $K = \C(S)$ of a smooth algebraic curve $S/\C$. Here is a typical corollary of our main theorem.

\medskip
Let $\pi : {\bf B} \rightarrow S$ be a semi-abelian scheme
of constant toric rank,
 with generic fiber $B/K$, Lie algebra $LB/K$, let $A/K$ be the maximal abelian variety occuring as a quotient of $B$, and let $(A_0, \tau)$ be the $K/\C$ trace of $A/K$, where we assume for simplicity that $A_{0}$ is included in $A$.  By base change to a finite cover of $S$, we may assume that it also is the $K^{alg}/\C$-trace of $A/K^{alg}$, 
and we will just call it the $\C$-trace, or the constant part, of $A$.  Let further $B_0$ be the pull-back of  
$B$  to $A_0$. This is a semi-abelian variety, which is defined over $K$ and does not in general descend to $\C$; we call $B_0$  the semi-constant part of $B$.  Finally, let $\tilde {\bf B}$ be the universal vectorial extension of $\bf B$, with generic fiber $\tilde B/K$. 
Consider the exponential morphisms of analytic sheaves attached to the group schemes    ${  {\bf B}}^{an}, {\tilde {\bf B}}^{an}$ over $S^{an}$ :
$$exp_{  {\bf B}}:  L{  {\bf B}}^{an} \longrightarrow { {\bf B}}^{an} ~, ~exp_{\tilde {\bf B}}:  L{\tilde {\bf B}}^{an} \longrightarrow {\tilde {\bf B}}^{an}.$$
Let further $x$ be a point in $LB(K)$. Restricting $S$ if necessary, we may assume that $x$ extends to a section ${\bf x}\in L{\bf B}(S)$, and lift it to a section $\tilde {\bf x}$ of $L\tilde {\bf B}(S)$, for which we set $ {\bf y} = exp_{{ {\bf B}}}(  {\bf x}) \in { {\bf B}}^{an}(S^{an}), \tilde {\bf y} = exp_{{\tilde {\bf B}}}(\tilde {\bf x}) \in {\tilde {\bf B}}^{an}(S^{an})$. Abbreviating these analytic sections as $y = exp_B(x), \tilde y= exp_{\tilde B}(\tilde x)$, 
and making use of the same notations over all covers of $S$, we will prove:

\begin{Corollary} Let $B/K^{alg}$ be a semi-abelian variety, and  let $x$ be a point in $LB(K^{alg})$.  Assume that

\medskip
\noindent
$({\bf HB}_0)$: the semi-constant part $B_0$ of $B$ is defined over $\C$;

\medskip
\noindent
$({\bf HX/B})_{K^{alg}}$:  for any proper semi-abelian subvariety  $H$ of $B$, $x \notin LH(K^{alg}) + LB_0(\C)$.

\medskip
\noindent
Let further $\tilde x \in L\tilde B(K^{alg})$ be a lift of $x$. Then, the fields of definition  of the points $y = exp_B(x), \tilde y = exp_{\tilde B}(\tilde x)$  satisfy:
$$ tr.deg.(K^{alg}(\tilde y)/K^{alg}) = dim(\tilde B),$$
and in particular :  $ tr.deg.(K^{alg}(  y)/K^{alg}) = dim( B).$ 
\end{Corollary}

\noindent
(We recall that by rigidity, all abelian subvarieties $H$ of $B$ are defined over $K^{alg}$.) 

\medskip
Consider for instance the case where $B = A$ is an abelian variety, so that $({\bf HB}_0)$ is automatically satisfied. When $A $  is a power of an elliptic curve,  Corollary 1.1 says that if $\wp$ is an elliptic function with a non constant invariant $j \in \C(z)$ and zeta function $\zeta$, and if $x_1(z), ..., x_n(z)$ are $\Z$-linearly independent algebraic functions, then the $2n$ analytic functions defined on some open domain of $\C$ by $\wp(x_1(z)), ..., \wp (x_n(z)), \zeta(x_1(z)), ..., \zeta(x_n(z))$ are algebraically independent over $\C(z)$. In the opposite case where $j \in \C$, this holds only if $x_1, ..., x_n$ are linearly independent, modulo $\C$, over the ring of multipliers of $\wp$.  When, as in the latter case, the full abelian variety 
 is constant ($A = A_0$), such results follow from Ax's work on the Schanuel conjecture, cf \cite{Ax}, Theorem 3, 
\cite{Kirby}, \cite{Bertrand} - and in the elliptic case, \cite{Brownawell-Kubota}.

\medskip
\begin{Remark} A similar result holds for the exponential morphism  of the {\it formal group} $\hat {\tilde B}$  of $\tilde B$ at the origin, cf. Appendix, Remark 7.6.  
Namely, assuming that the section  ${\tilde {\bf x}}\in L{ {\tilde {\bf B}}}(S)$ vanishes at a given point $s_0$ of $S$, and denoting by $\hat S_{s_0}$ the formal completion of $s_0$ in $S$, we can deduce from our methods  a direct proof that ${\hat{\tilde {\bf y}}} = exp_{\hat{\tilde {\bf B}}}({\hat{\tilde {\bf x}}}) \in {{\hat{\tilde {\bf B}}}} ({\hat S}_{s_0})$ has transcendence degree $dim(\tilde B)$ over $K$. Note that  this formal setting is the framework used by  Ax in the constant case, cf. \cite{Ax}, Theorem 3.
In an analogous way, one can replace in this corollary $\C$ by the $p$-adic field $\C_p$, with a suitable convergence condition on the point  $x$.

\end{Remark}

\medskip
As in Ax's initial work on the toric case, Corollary 1.1 will follow from consideration of the differential relations or equations satisfied by $(x,y)$ such that $exp(x) = y$. But
contrary to Ax's setting,  the ambient semiabelian variety $B$ will here not be constant, and we will in general only be able to find the appropriate differential equations on the universal vectorial extension $\tilde B$ of $B$ (and its Lie algebra). Our main theorem, on solutions $\tilde y$ of $\partial\ell n_{\tilde B}(\tilde y) =
\partial_{L\tilde B}(\tilde x)$, is Theorem 1.4 below. The inductive nature of the proof (passing to quotients)  will force us into the more general category of {\em almost semiabelian $D$-groups} and our main technical result is Theorem 1.3 below.  More precisely, and referring to the following 
sections of the paper for the underlying notions, we will study the following problem.

\medskip
Let 
now $K
= K^{alg}$ be the algebraic closure of a function field in one variable over $\C$, let $\partial$ be a non-trivial derivation on $K$, let $\mathcal U$ be some universal differential extension of $K$, and 
 let $K^{diff}$ be a differential closure of $K$ in $\mathcal U$. 
 More will be said about $K^{diff}$ and ${\cal U}$ at the beginning of section 2. For now we just mention that $K^{diff}$ and $K$ 
have the same field of constants, namely $\C$.
 
Let $G/K$ be a connected {\it commutative   algebraic $D$-group} (cf.  \cite {Buium}, \cite{Pillay}  and Section 2 below), i.e. $G$ is a connected commutative algebraic group over $K$, equipped with an extension to ${\mathcal O}_{G}$ of the derivation $\partial$, which respects the group structure of $G$. Denoting Lie algebras by  $L$, we write:
$$\partial \ell n_{G} : G \rightarrow LG$$
for the corresponding {\it logarithmic derivative} on $G$ (cf. Section 2). This is a
``first order
differential algebraic"
homomorphism which takes $G(\mathcal U)$ onto $LG(\mathcal U)$, and likewise takes $G(K^{diff})$ onto $LG(K^{diff})$, but will be far from surjective at the level of $K$-points. The kernel of $\partial\ell n_{G}$, denoted   $G^\partial$ when the $D$-group structure on $G$ is assumed, is a ``differential algebraic group" defined over $K$, and for any differential extension field $K'$ of $K$ (including $K^{diff}$ and ${\mathcal U}$), we can speak
of the group $G^{\partial}(K')$ of its $K'$-points.
We will write:
$$\partial_{LG} : LG \rightarrow LG$$
for the canonical connection, contracted with $\partial$, which $\partial \ell n_{G}$  induces on $LG$, and which
we can again view as a differential algebraic endomorphism of $LG$, surjective at the level of ${\mathcal U}$-points. This is discussed in detail in the Appendix of this paper; for instance, when $G$ is the universal vectorial extension of an abelian variety $A$, $\partial_{LG}$ coincides with the dual of the standard Gauss-Manin connection on $H^1_{dR}(A/K)$. Again we  write  $(LG)^\partial$ for the kernel of $\partial_{LG}$, namely
the space of vectors horizontal for the connection $\partial_{LG}$.

\bigskip

We consider  the differential relation 
$$ \partial \ell n_{G} (y) = \partial_{LG} (x)  \quad (*) ,$$
where $(x, y) \in  (LG Ê\times G )(\mathcal U)$, and proceed to compute the transcendence degree of $K(y)$ over $K$ under the assumption that $x$ is $K$-rational. As will 
later become apparent, this is the natural algebraic description of  the ``Lindemann-Weierstrass case" of the Schanuel conjecture:  
with some abuse of notations, the rough idea is that whenever $exp_G(x)$ is well-defined, we have
$$\partial \ell n_G (exp_G(x)) = \partial_{LG}(x)$$
(see the proof of Corollary 
1.1 at the end of this introduction, and Section G of the Appendix), so that up to addition by elements $x_0 \in (LG)^{\partial}({\cal U}), y_0 \in G^{\partial}({\cal U})  $,
$$ \partial \ell n_{G} (y) = \partial_{LG} (x) \Leftrightarrow  y-y_0 = exp_G(x-x_0).$$
In the present paper, we do not discuss 
the case where  $y$ is $K$-rational and 
$x$ is the unknown, nor the general case where both $x$ and $y$ are unknown. These would respectively correspond to the functional analogues of the Grothendieck and 
Schanuel-Andr\'e
 conjectures, cf. \cite{Bertrand}. 

\medskip
We make the necessary assumption that our algebraic $D$-group $G$  admits no non-
zero  vectorial quotient, and thereby restrict to {\it almost semi-abelian $D$-groups}, as defined in Section 3 :  these are the quotients 
$$ G = \tilde B/U$$
by a vectorial $D$-subgroup  $U$ of the universal extension $\tilde B$ of a semi-abelian variety $B$ 
(the latter being endowed with its unique $D$-group structure). 
Furthermore, we will assume that $B$ satisfies Hypothesis $({\bf HB})_0$ of Corollary 1.1, and will rename this hypothesis as  

\medskip
\noindent
$({\bf HG})_0$: {\it  the semi-constant part $B_0$ of the maximal semi-abelian quotient $B$ of $G$  is actually constant.}

\medskip
\noindent
Let $B_{(0)}$ be the constant part of $B$, i.e. the maximal 
semiabelian subvariety of $B$ isomorphic to one defined over $\C$. In general, $B_{(0)}$ is contained in $B_0$, and hypothesis $({\bf HG})_0$ means that they coincide.  
In particular, it is automatically satisfied if the maximal abelian quotient $A$ of $B$ is traceless (in which case $B_{(0)} = B_0$ is the toric part of $B$),  or if  $B$ is the product of the abelian variety $A$ by a torus.  Our main 
technical result can then be stated as follows.

 \begin{Theorem} : let  $G$  be an almost  semi-abelian $D$-group,
 defined over $K$,  which  satisfies Hypothesis $({\bf HG})_0$,
and let $x$ be a point  in  $LG(K)$. Assume that
 
 \medskip
\noindent
$({\bf HX})_K$ : $x\notin LH(K) + (LG)^{\partial}(K)$ for any proper algebraic subgroup $H/K$
of $G$.

\medskip
\noindent
 Let $y \in G({\mathcal U})$ be a solution of the differential equation $ \partial \ell n_{G} (y) = \partial_{LG} (x)$. Then,   
 $$tr.deg. (K(y)/K) = dim(G).$$

\end{Theorem}

We will give 
(see \S 5.3) an example  showing that the $({\bf HG})_0$ hypothesis in Theorem 1.3 cannot be dropped. 
On the other hand, Theorem 5.2 shows that  
 hypothesis $({\bf HG})_0$ {\em can} be dropped, but at the expense of strengthening the hypothesis
$({\bf HX})_{K}$ to 
$({\bf HX})$: $x\notin LH({\cal U}) + LG^{\partial}({\mathcal U})$ for any proper algebraic subgroup $H/K$ of $G$. However we tend to prefer the 
$({\bf HX})_{K}$ hypothesis
because, as witnessed by Corollary 1.1, the results it yields are closer in spirit to the number theoretic case.

\medskip
\noindent
The proofs of Theorem 1.3 and these other results will be given in section 5 of the paper. In \cite{Bertrand} it was suggested that differential Galois theory in the most general form, may be useful in this function field Lindemann-Weierstrass context. Under an additional assumption ``$K$-largeness" on the algebraic $D$-group $G$ (which actually implies $({\bf HG})_{0}$), such a Galois-theoretic proof of 1.3 is in fact possible, and is given in section 6. 

\medskip
When $G = \tilde B$ is the universal extension of a semi-abelian variety $B$ as in Corollary 1.1,  Theorem 1.3 yields 
our main result: 
\begin{Theorem} :  Let $B/K$ be a semi-abelian variety, and  let $x$ be a point in $LB(K)$.  Assume that

\medskip
\noindent
$({\bf HB}_0)$: the semi-constant part $B_0$ of $B$ is defined over $\C$;

\medskip
\noindent
$({\bf HX/B})_{K}$:  for any proper semi-abelian subvariety  $H$ of $B$, $x \notin LH(K) + LB_0(\C)$.

\medskip
\noindent
Let further $\tilde x \in L\tilde B(K)$ be a lift of $x$, and let $\tilde y \in \tilde B({\mathcal U})$ be a solution of the differential equation $ \partial \ell n_{\tilde B} (\tilde y) = \partial_{L\tilde B} (\tilde x)$. Then,   
 $$tr.deg. (K(\tilde y)/K) = dim(\tilde B).$$
 \end{Theorem}

Since the the almost semi-abelian $D$-groups of Theorem 1.3 are quotients of such $D$-groups $\tilde B$, it is conversely  clear that Theorem 1.4 implies Theorem 1.3 (details are given below).

\subsection{Organisation  of the proofs}

The proof of Theorem 1.3 (and variants) in section 5 has a number of ingredients and draws on several sources. On the one hand there is differential algebraic geometry (and model theory) which provides the notions of algebraic $D$-group, logarithmic derivative, etc.,  in which our results are phrased
in the main body of the paper,
 as well as the powerful ``socle theorem" which is closely related to the  function field Mordell-Lang conjecture in characteristic $0$ , and in the present context facilitates an inductive proof. Then there are results originating with Ax dealing with the case of $G$ defined over $\C$. Finally, from algebraic geometry we make use of the Manin-Coleman-Chai theorem of the kernel, as well as the 
 Griffiths-Schmid-Deligne theorem of the fixed part and Deligne's semisimplicity theorem. To be able to draw on these various sources we need to know at least the
{\em compatibility} of the different languages and constructions. Among the issues is the relation between the logarithmic derivative on the universal extension $\tilde A$ of an abelian variety $A$ defined over $K$, and the Gauss-Manin connection on $H^{1}_{dR}(A/K) =$  dual of $L(\tilde A)$.  
So our rather extensive appendix, {\it Exponentials on Algebraic $D$-groups},  is devoted to clarifying some of these issues,
although they are probably well-known. A discussion of  the ``theorem of the kernel" also appears there.

\medskip
\noindent 
In sections 2 and 3 we introduce and discuss algebraic $D$-groups, 
differential algebraic groups, $D$-modules, and 
almost semiabelian $D$-groups, as well as logarithmic derivatives, in the context of Kolchin's differential algebraic geometry. Section 3 contains a few new observations.
Section 4 presents the main tools (including the ``socle theorem")  and special cases
 which will be used in the proof of the main theorem (1.3 above). 
As already said, the proof of the main theorem, plus a variant as well
as a (counter)example, are given in section 5, while  Section 6 provides
another proof in the special case where the algebraic $D$-group $G$ is
$K$-large, based on differential Galois theory in place of the socle
theorem.  
As for the Appendix, its main results have been gathered in a Conclusion before its section J.

\medskip
We conclude this introduction by discussing the mutual relations between Theorem 1.3, Theorem 1.4 and Corollary 1.1, in particular showing how to deduce
Corollary 1.1 from Theorem 1.3. 

\medskip
\noindent
{\it Theorem 1.3 $\Leftrightarrow$ Theorem 1.4}

\medskip

Since in the projection $\tilde B \rightarrow \tilde B/U = G$, transcendence degrees decrease at most by the dimension $U$ of the kernel, Theorem 1.3 on a general $G$ is equivalent to its special case for $\tilde B$.  Since  Hypothesis $({\bf HB}_0)$ is just a renaming of $({\bf HG}_0)$, it remains, given a point $x \in LB(K)$ and a lift $\tilde x$ of $x$ to $L\tilde B(K)$,  to check that $x$ satisfies hypothesis $({\bf HX/B})_{K}$ of Theorem 1.4 with respect to $LB$ if and only if $\tilde x$  satisfies the corresponding hypothesis, say $({\bf H\tilde X})_K$, of Theorem 1.3 with respect to $L\tilde B$.  Now, we will show in 
Lemma 4.2(i)
(see Corollary H.5.ii of the appendix in the case $B$ is an abelian variety) that under  $({\bf HB}_0)$, the universal vectorial extension $\tilde B_0 \subset \tilde B$ of $B_0$, which is then clearly  defined over $\C$, satisfies:
$$(L\tilde B)^{\partial}(K) = L\tilde B_0(\C).$$
 Furthermore, $\tilde B$ is an essential extension of $B$, that is to say: any proper algebraic subgroup  of $\tilde B$ projects onto a proper algebraic subgroup of $B$.

So, 
suppose that there exists a proper algebraic subgroup $H'/K$ of $\tilde B$ with projection $H$ on $B$,  such that $\tilde x$ lies in $LH'(K) + L{\tilde B_0}(\C)$.  Projecting to $B$, we deduce from $({\bf HX/B})_{K}$ that $H$ fills up $B$, and  $H'$ cannot be proper. Conversely,   suppose that $x$ lies in $LH(K) + LB_0(\C)$ for some proper semi-abelian subvariety $H$ of $B$, with inverse image $H'$ in $\tilde B$;  then $H'$ is a proper algebraic subgroup of $\tilde B$ over $K$ satisfying $\tilde x \in LH'(K) + (L\tilde B)^{\partial}(K)$, which contradicts $({\bf H\tilde X})_K$. 

\medskip
\noindent
{\it Theorem 1.4 $\Rightarrow$ Corollary 1.1}

\medskip

We go back  to the notations before the statement of this
 corollary, but now denote by  $K = \C(S)^{alg}$
 the algebraic closure of the field $\C(S)$.

 So, let $\tilde {\bf y} := exp_{\tilde {\bf B}}(\tilde {\bf x})  \in \tilde {\bf B}(S^{an})$ be the image of $\tilde {\bf x}$ under the exponential morphism of the group scheme $\tilde {\bf B}^{an}/S^{an}$. 
  By one of the compatibilities proven in the Appendix, cf. Proposition G.1, its logarithmic derivative
   $\partial \ell n_{\tilde B}$, extended to $\tilde {\bf B}/S$, satisfies
$$\partial \ell n_{\tilde {\bf B}} \tilde {\bf y}= \partial \ell n_{\tilde {\bf B}} (exp_{\tilde {\bf B}} \tilde {\bf x}) = \partial_{L\tilde {\bf B}}  \tilde {\bf x}.$$
Viewing the field of meromorphic functions on $S^{an}$ as a subfield of $\cal U$,  we have thus constructed a solution $\tilde y \in \tilde B({\cal U})$ of the differential equation    $\partial \ell n_{\tilde B}  \tilde y = \partial_{L\tilde B} \tilde x$,  with $\tilde x \in \tilde B(K)$. Theorem 1.4 then implies that the transcendence degree of $K(\tilde {\bf y}) = K(\tilde y)$ over $K$ is equal to $dim(\tilde B)$. As for the last statement of the corollary, it again follows from the fact that in the  projection from $\tilde B$ to $B$, transcendence degrees can at most decrease by the dimension of the kernel.

Notice that contrary to $\tilde B$,  the semi-abelian variety $B$ admits in general no $D$-group structure, so that  the relation $y = exp_B(x)$ cannot be expressed directly on $B$ in 
differential algebraic terms.  In other words,  even for the study of the point $y$, 
lifting to $\tilde B$ is forced onto us in order to allow for the techniques of proofs of Theorem 1.3.

\medskip
\begin{Remark} We are aware that statements like Corollary 1.1 can  often be proved by purely analytic means, 
by using the order of growth of $exp_{\tilde {\bf B}}(\tilde {\bf x})$ at the essential singularities given by the poles of the rational section $\tilde x \in L\tilde B(K)$. But as mentioned in  Remark 1.2, the scope of our methods is broader. Corollary 1.1 is merely an illustration, while the true object of study of this article is the differential relation $(*)$ in differential fields.
 
\end{Remark}

\vspace{5mm}
\noindent
{\em Acknowledgement.}  Both authors would like to thank the European model theory network MODNET. The 2006 meeting in Antalya was the starting point for this research project and subsequent MODNET meetings in Luminy (2007) and La Roche (2008) gave us the opportunity to meet and continue the collaboration in person.

\section{Differential algebraic preliminaries}

The context here is the differential algebra or differential algebraic geometry of Ritt and Kolchin, as developed in Kolchin's books  \cite{Kolchin1}, \cite{Kolchin2}.  We refer to Buium's books \cite{Buium} (especially Chapter 5) and \cite{Buium2}, to the second author's paper \cite{Pillay} and to Malgrange's monograph \cite{Malgrange}  for the specific notions needed for the current paper, although we will give a brief account in this section.

We fix a universal differential field $({\cal U},\partial)$ of characteristic $0$ in which all differential fields we discuss are assumed to embed. We denote by ${\cal C}$  the field of constants of ${\cal U}$, 
by $K$   an 
algebraically closed
differential field (differential
subfield of ${\cal U}$),   by $C_{K}$ its field of constants and by $K^{diff}$  a differential closure of $K$. Of course in our main applications $K$ will have transcendence degree $1$ over $C_{K} = \C$.  
Here is a quick description of ${\cal U}$ and differential closures: A differential field $L$ is said to be differentially closed if any finite system of differential polynomial equations over $L$ in unknowns $x_{1},..,x_{n}$ which has a common solution in some differential field extending $L$ already has a solution in $L$.
Differentially closed fields exist. 
More precisely,
 fix an uncountable cardinal $\kappa$. Then ${\cal U}$ 
will be  a differentially closed field of cardinality $\kappa$ with the following property: whenever $L_{1}< L_{2}$ are differential fields of cardinality $< \kappa$ and $f:L_{1}\to {\cal U}$ is an embedding (of differential fields) then $f$ extends to $g:L_{2} \to {\cal U}$. In our context, $\kappa$ is assumed to be strictly greater than the cardinality of our base differential field $K$. A {\em differential closure} of $K$ is a differentially closed field extending $K$ with the property that it embeds over $K$ into any differentially closed field containing $K$. Again a differential closure of $K$ exists and is unique up to isomorphism over $K$. So $K^{diff}$ denotes a differential closure of $K$ inside ${\cal U}$.

\medskip
\noindent
The geometric objects of Kolchin's theory are ``differential algebraic varieties", which are given locally as 
common solution sets in ${\cal U}^{n}$ of finite systems of differential polynomial equations. :
They form a category, whose morphisms are easily defined and will be termed  ``differential".
Furthermore, we say that  a differential algebraic variety $X$ is defined over $K$ if its defining equations have coefficients in $K$. One can then view $X$ as a functor from the category of differential $K$-algebras to sets.  
 Likewise one has the notion of a {\it differential algebraic group}.
 However the differential algebraic groups we consider will all 
be given explicitly as subgroups of algebraic groups.
If $X$ is a differential algebraic variety  (in particular if it is an algebraic variety) defined over a subfield of ${\cal U}$,
 we will 
 often identify $X$ with its set $X({\cal U})$ of ${\cal U}$-points.

\subsection{The twisted tangent bundle}
If $X\subseteq {\cal U}^{n}$ is an affine algebraic variety, and $a\in X$
, we can apply $\partial$ to the coordinates of $x$ to obtain
a point $\partial(x)\in {\cal U}^{n}$. 
This depends of course on the chosen embedding $X\subseteq {\cal U}^{n}$, but it can be viewed in an intrinsic way 
as a (differential rational rather than rational) section of a certain {\em twisted} tangent bundle $T_{\partial}(X)$ of $X$, which we now describe. 
We assume $X$ smooth, (geometrically) irreducible, and defined over $K$. If 
$X\subseteq {\cal U}^{n}$ is affine, then $T_{\partial}(X) = \{(a,u)\in {\cal U}^{2n}: a\in X$ and 
$\sum (\partial f/\partial x_{i})(a)u_{i} + f^{\partial} = 0$ for $f$ ranging over generators of $I(X)\}$. Here 
$f^{\partial}$ is obtained from $f$ by applying $\partial$ to the coefficients of $f$. For arbitrary $X$, take a covering by open affines $U_{j}$ and piece together the $T_{\partial}(U_{i})$ using the transition maps in the obvious way, to obtain $T_{\partial}(X)$. So if $X$ is defined over $C_{K}$, then $T_{\partial}(X)$ coincides with the tangent bundle $T(X)$ of $X$. In general the definition shows $T_{\partial}(X)$ to be a torsor (over $X$) under
 $T(X)$.  For $X$ affine one sees from the Leibniz rule that if $a\in X$ then $(a,\partial(a))\in T_{\partial}(X)$.
In the general case this makes sense too, and with abuse of notation we call $\partial: X\to T_{\partial}(X)$ the corresponding (differential regular)
section. If $s: X\to T_{\partial}(X)$ is a regular section, defined over $K$, then we obtain a (finite-dimensional)
 {\em differential algebraic variety} $(X,s)^{\partial} = \{a\in X: 
\partial(a) = s(a)\}$, defined over $K$. Finite-dimensionality means the following: Suppose $X\subseteq {\cal U}^{n}$ is an affine differential algebraic variety defined over $K$. We call $X$ finite dimensional if there is 
a positive integer $m$ such that for any point $a\in X$, the differential field  $K\langle a\rangle$ generated by $K$ and $a$ has transcendence degree at most $m$ over $K$.
So if $X$ is an {\em algebraic} variety (i.e. with no addditional differential equations)
and not a point, 
 then $X$ is infinite-dimensional. On the other hand, clearly our $(X,s)^{\partial}$ above is finite-dimensional, as for $a\in (X,s)^{\partial}$, $K\langle a \rangle = K(a)$.

\medskip
In sections A and B of the appendix,  a geometric  account of these and the next notions is given when the base $K$ is replaced by a curve $S$ over $\C$, and  $\partial|K$ 
by
a vector field on $S$.  In fact (like $T$)
 $T_{\partial}$ is a functor:
if $\phi:X\to Y$ is a morphism over $K$ (of smooth irreducible varieties), then $T_{\partial}(\phi): T_{\partial}(X) \to T_{\partial}(Y)$ is given in local coordinates by $T\phi + \phi^{\partial}$.

\medskip
\noindent

\subsection{Algebraic $D$-groups
and logarithmic derivatives}
If $G$ is now a connected algebraic group defined over $K$, then because $T_{\partial}$ is a functor, $T_{\partial}(G)$ 
(also denoted $T_\partial G$) 
can be naturally equipped with the structure of an algebraic group over $K$, and the canonical projection $\pi:T_{\partial}(G)\to G$ is a homomorphism of algebraic groups, whose kernel is canonically isomorphic to the Lie algebra $LG$ of $G$. Indeed, this kernel is 
$(T_{\partial}G)_{e}$, the fibre of $T_{\partial}G$ above the identity $e$ of $G$, which we have seen is a principal homogeneous space for $T(G)_{e} = LG$ over $K$. 
Now the identity element of $T_{\partial}G$, that is the $K$-rational point $(e,\partial(e))$, gives rise to an identification of $T_{\partial}G_{e}$ and $LG$   (See section 2 of \cite{Pillay}.)

By an {\it algebraic $D$-group}
(or just: $D$-group)
 over $K$, we mean a pair $(G,s)$ where $G$ is an algebraic group over $K$ and $s:G\to T_{\partial}(G)$ is a regular section defined over $K$ which is also a homomorphism. Algebraic $D$-groups are objects of algebraic geometry. Giving $G$ the structure $s$ of an algebraic $D$-group over $K$ is {\em equivalent} to extending the derivation $\partial|K$ to a derivation of the structure sheaf of $G$ over $K$ which respects the group operation, and this is how they were first defined by Buium \cite{Buium}.

We will restrict our attention to commutative algebraic groups $G$. In this case $T_{\partial}(G)$ is also commutative, and we will write its group law using additive notation. One can check that $\partial: G\to T_{\partial}(G)$ is also a homomorphism, so if $(G,s)$ is an algebraic $D$-group, then $\partial - s$  (where $-$ is meant in the sense of the group $T_{\partial}(G)$) is a {\em differential regular} homomorphism from $G$ to $LG$, which we call the  {\it logarithmic derivative} associated to $(G,s)$ and which is written as $\partial \ell n_{(G,s)}$ or
$\partial \ell n_{G}$ when $s$ is understood. Consistent with earlier notation we write $(G,s)^{\partial}$, or just $G^\partial$,
 for the kernel of this logarithmic derivative. $(G,s)^{\partial}$ {\em is} a (finite-dimensional) differential algebraic group. 
If $(G,s)$ is defined over $K$, then for any differential field $F$ containing $K$, we have  $\partial\ell n_{(G,s)}: G(F) \to LG(F)$. For $F = {\cal U}$ or $K^{diff}$, this map is surjective.

\medskip

Two special cases of an algebraic $D$-group should be familiar. First when $G$ is defined over $C_{K}$ and $s:G\to T(G)$ is the $0$-section. Then the corresponding logarithmic derivative reduces to the classical one of Kolchin, so for example if $G$ is 
a torus, it is just $\partial(g)\cdot g^{-1}$. Moreover $G^{\partial}
(\cal U)$
 is just $G({\cal C})$.

The second is when $G = \G_{a}^{n}$. Then $LG = G$ and $T(G) = G \times G$. A regular section $s$ defined over $K$ is of the form $x \mapsto (x,Ax)$ where $A$ is an $n\times n$ matrix over $K$ and we write $x\in G$ as a column vector.
Hence the corresponding logarithmic derivative is 
$x \mapsto \partial(x) - Ax$ from $G$ to $G$. This is precisely
a $D$-module structure on ${\cal U}^{n}$ or on $K^{n}$ if we restrict to $K$-points, in other words a linear differential system. 
Indeed,  a $D$-module defined over $K$ is by definition a $K$-vector space $V$  (that is, a commutative unipotent group over $K$ as we are in characteristic $0$) together with an additive homomorphism $D_{V}: V \to V$ satisfying, with $\partial$, the Leibniz rule. 
The algebraic group
 $V$ is isomorphic to some $\G_{a}^{n}$ over $K$, so the $D$-module structure is as above the logarithmic derivative for some algebraic $D$-group structure on $\G_{a}^{n}$.

\subsection{$D$-modules and the connection $\partial_{LG}$ on $LG$}
We discuss here the $D$-module structure induced on the Lie algebra of an algebraic $D$-group $(G,s)$, which yields our map $\partial_{LG}$ referred to in the introduction. This is a rather delicate aspect of the paper, in terms of compatibilities, and much of the Appendix is devoted to it. Here we give an ``algebraic" definition based on Section 5 of \cite{Pillay-Ziegler}
 and point out some equivalences. 
First suppose $X$ to be a smooth irreducible affine variety over $K$ with coordinates $x = (x_{1},..,x_{n})$.
Recall from 
subsection 2.1 
that the tangent bundle of $X$ is defined by equations for $X$ together with $df\cdot y = 0$ for $f$ ranging over generators of the ideal of $X$ (where $df\cdot y = \sum_{i=1,..,n} (\partial f(x)/\partial x_{i})y_{i} =  0$ represents the (vertical) differential 
$d_{G/K}f$ of $f$).  On the other hand the twisted tangent bundle of $X$ is defined by equations defining $X$ together with  $df(x)\cdot y + f^{\partial}(x)$ for $f$ generating the ideal of $X$.
So we choose ($x,u)$ as coordinates of $T(X)$ and $(x,y)$ as coordinates of $T_{\partial}(X)$. 
Likewise let $(x,u,y,v)$ be coordinates for $T(T_{\partial}X)$ and  $(x,y,u,v)$ coordinates for $T_{\partial}(T(X))$.
A straightforward computation using the formulas above yields that the map taking $(x,y,u,v)$ to $(x,u,y,v)$ gives an isomorphism between $T(T_{\partial}X)$ and $T_{\partial}(T(X)))$ defined over $K$. By gluing, this extends to
arbitrary smooth varieties over $K$. We obtain (see \cite{Pillay-Ziegler}, Lemma 5.1):
\begin{Lemma} Let $G$ be a commutative algebraic group over $K$. Then there is a functorial isomorphism between the commutative algebraic groups $T(T_{\partial}G)$ and $T_{\partial}(T(G))$ (over the natural projections to $T_{\partial}(G))$, which yields also a functorial isomorphism between $L(T_{\partial}G)$ and $T_{\partial}(LG)$.
\end{Lemma}

Suppose now that $(G,s)$ is a commutative algebraic $D$-group over $K$. Then from $s:G\to T_{\partial}G$ we obtain (differentiating at the identity
and setting $d_{G/K, e} = L$), a homomorphism $Ls:LG\to L(T_{\partial}G)$, which via the identifications of Lemma 2.1  gives a (regular) homomorphic section $Ls: LG \to T_{\partial}(LG)$, that is a $D$-group structure on $LG$, all defined over $K$. The corresponding logarithmic derivative on $LG$  is
$\partial - Ls:  LG \to LLG = LG$, namely $\partial \ell n_{(LG,Ls)}$, or $\partial \ell n_{LG}$ when $Ls$ is assumed. This is the $D$-module structure on $LG$ that we are interested in, and we will denote it
$\partial_{LG}$ for simplicity. 
$(LG)^{\partial}$ will denote the kernel of $\partial_{LG}$. For a geometric account, see sections C and D of the appendix.

\begin{Remark} If $G$ has dimension $n$, is defined over $C_{K}$ and $s = 0$, then $(LG, Ls)$ is isomorphic to $(\G_{a}^{n}, 0)$ and 
 $\partial_{LG} = \partial$ on $\G_{a}^{n}$.
\end{Remark}

We point out two alternative algebraic descriptions of $\partial_{LG}$. The first is as in \cite{Pillay-Ziegler}, 
 Lemma 3.7 (iii) and the paragraph following it:  we know that $s:G\to T_{\partial}G$ gives a derivation, which we still call $s$, of the structure sheaf of $G$ over $K$ extending $\partial|K$. Now, $s$ acts on the local ring of $G$ at the identity, and in fact preserves the maximal ideal ${\cal M}$. So $s$ induces a $D$-module structure on the cotangent space ${\cal M}/{\cal M}^{2}$ at identity, and hence the dual connection on $LG$. This can be checked to coincide with $\partial_{LG}$, as mentioned  in \cite{Pillay-Ziegler}, Lemma 5.1 (ii).

Another description is via differentials of differential regular functions, as defined by Kolchin \cite{Kolchin2}, 
Section 2 of Chapter 8. We have the
differential regular homomorphism $\partial \ell n_{G}:G\to LG$, which has a differential $L\partial \ell n_{G}$ at the identity, a differential regular homomorphism from $LG$ to $LG$ which again can be shown to coincide with $\partial_{LG}$.  

\subsection{Algebraic $D$-groups and differential algebraic groups}
We recall some results and facts from \cite{Kowalski-Pillay}.
Suppose that $(G,s_{G})$, $(H,s_{H})$ are connected commutative algebraic $D$-groups defined over $K$, and $h:G\to H$ is a homomorphism of algebraic groups defined over $K$. We will say that $h$ is a homomorphism of algebraic $D$-groups, if $T_{\partial}h\circ s_{G} = s_{H}\circ h$. 

If $H$ is an algebraic subgroup of the algebraic $D$-group $(G,s)$, we say that $H$ is a $D$-subgroup if 
$s|H: H \to T_{\partial}H \subseteq T_{\partial}G$. (Likewise for ``$D$-subvariety".)

From 2.7 of \cite{Kowalski-Pillay} we obtain (suppressing mention of $s$ sometimes)
the easy
\begin{Fact} (i) Let $h:G\to H$ be a homomorphism of algebraic $D$-groups. Then $Ker(h)$ is an algebraic $D$-subgroup of $G$.
\newline
(ii) Conversely if $H$ is an algebraic $D$-subgroup of $G$, then $G/H$ can be equipped with the structure of a $D$-group such that the quotient map is a $D$-homomorphism.
\end{Fact}

From Fact 2.3 of \cite{Kowalski-Pillay} we have the deeper 
\begin{Fact} (i) 
If $G$ is an algebraic $D$-group then $G^{\partial}$ is Zariski-dense in $G$. Consequently, if $(G,s_{G})$, $(H,s_{H})$ are connected algebraic $D$-groups, and $h:G\to H$ is a homomorphism of algebraic groups, then $h$ is a $D$-homomorphism if and only if $h(G^{\partial}) \subseteq H^{\partial}$.
\newline
(ii) Let $G$ be a connected algebraic $D$-group. Then we have a bijection between (connected) algebraic $D$
-subgroups of $G$ and (connected) differential algebraic subgroups of $G^{\partial}$: namely for $H$ a $D$
-subgroup of $G$, $G^{\partial}\cap H = H^{\partial}$ is a differential algebraic subgroup of $G^{\partial}$ whose Zariski 
closure is $H$. And if $\cal H$ is a differential algebraic subgroup of $G^{\partial}$, then the Zariski closure  of $\cal H$ 
in $G$ is a $D$-subgroup  whose intersection with $G^{\partial}$ is precisely $\cal H$.
\end{Fact} 

\bigskip
We will also need the following easy fact 
(a converse statement for embeddings is given in Corollary G.5 of the appendix):
\begin{Lemma} Suppose that $f:G\to H$ is a homomorphism of algebraic $D$-groups, then $Lf:LG \to LH$ is a homomorphism of $D$-modules (where $LG$, $LH$ are equipped with $D$-module structures as in section 2.3).
\end{Lemma}

\medskip
\noindent
{\bf Note.} - Throughout this section 2.4, we can also specify fields of definition where appropriate.  For example in Fact 2.3 (i) if $h$ is defined over $K$ so is $Ker(h)$, and in Fact 2.3 (ii) if $G,H$ are defined over $K$ so is $G/H$.

\section{Almost semiabelian $D$-groups}

We recall that we work in characteristic $0$,
 in the context of a universal differential field ${\cal U}$ and with a small differential subfield $K$. 
We maintain our general assumption that $K$ is algebraically closed, but it is not always needed. Although much of what we say is implicit or explicit in the literature (such as \cite{Buium} and  \cite{Marker}) we may give proofs, for the convenience of the reader.

\subsection{Almost semiabelian varieties and $D$-groups}

Recall that an abelian variety $A$ has a ``universal vectorial extension",
which we denote by $\tilde A$ throughout the paper, and which admits the following characterization:
 there is an exact sequence
of commutative connected algebraic groups:
$$0 \mapright{} W_A \mapright{}
     \tilde A \mapright{\pi} A \mapright{}
      0 ~$$ 
 such that $W_{A}$ is unipotent (i.e. a vector group), and such that for any extension $f: G\to A$ of $A$ by a vector group there is a unique homomorphism $h:\tilde A \to G$ such that $\pi = f \circ h$.
\newline
Moreover, if $A$ is defined over $K$, so is $\tilde A$. In fact $W_{A}$ is the dual of $H^{1}(A,{\cal O}_{A})$ so $\tilde A$ has dimension $2dim(A)$.

\medskip
Likewise any semiabelian variety $B$ has a universal vectorial extension $\tilde B$ with the same universal property as above. In fact if $B$ is an extension of the abelian variety $A$ by the algebraic torus $T$, then $\tilde B$ is $B\times_{A}{\tilde A}$, which is an extension of $B$ by $W_A$. Again if $B$ is defined over $K$ so is $\tilde B$.

\medskip

\begin{Lemma} Let $G$ be a commutative connected algebraic group defined over $K$. The following conditions are equivalent.

i) $G$ has no nonzero homomorphisms to $\G_{a}$;

ii) there exists a semiabelian variety $B/K$ and   a unipotent subgroup $U/K$ of $\tilde B$ such that $G$ is isomorphic to   ${\tilde B}/U$.

iii) the group $Tor(G)$ of torsion points of $G$ is Zariski-dense in  $G$.

\end{Lemma}

\pf  
\newline
$i) \Rightarrow ii)$ :  by Chevalley's theorem, $G$ is an extension of its maximal semiabelian quotient $B$ by a unipotent group, say $V$.  Let $h :\tilde B \to G$ be given by the universal properties of $\tilde B$. So $U = Ker(h)$ is unipotent. Let $H = h(\tilde B)$. So $H$ projects onto $B$, and thus $G/H$ is unipotent (a subgroup of $V$). As we are assuming $G$ has no unipotent quotients, $G = H$, and $G = \tilde B/U$ as required.
\newline
  $ ii) \Rightarrow iii) $ :
  since $Tor( B)$ is Zariski dense in $B$, and since $\tilde B$ is an essential extension of $B$, $Tor(\tilde B)$ too is Zariski dense in $\tilde B$, and the same property is satisfied by the quotient $G$ of $\tilde B$.
\newline
$iii) \Rightarrow i)$ : suppose that $f$ is a surjective homomorphism from $G$ to $\G_{a}$. Then $Tor(G) < Ker(f)$.

\begin{Definition} We will call $G$ {\rm almost semi-abelian (asa)} if it satisfies the equivalent conditions of Lemma 3.1.
\end{Definition}
\begin{Remark}
Inspired by Brion's concept of anti-affine groups \cite{Brion}, we could alternatively say that $G$ is anti-additive.  
We refer the reader to the first part of Hypothesis $(\bf H)$ of
section I of the Appendix  for a Betti version of Condition (iii) (which we could have phrased in $\ell$-adic  terms). 
\end{Remark}

\vspace{5mm}
\noindent
For the rest of this section we consider $D$-structures on $asa$ groups. By an {\it almost semi-abelian $D$-group}, we mean an algebraic $D$-group $(G,s)$ such that $G$ is almost semi-abelian. 

\begin{Lemma} (i) If $G$ is $asa$, then $G$ has at most one structure of a $D$-group. Moreover if $G$ is defined over $K$, then so is the $D$-structure, if it exists.
\newline
(ii) If $B$ is semiabelian, then $\tilde B$ has a (unique) structure of $D$-group.
\newline
(iii) If $B$ is semiabelian and defined over $K$, then $B$ has the structure of a $D$-group if only if $B$ descends to $C_{K}$.
\end{Lemma}
\pf (i) If $s_{1},s_{2}$ were distinct rational homomorphic sections $G \to T_{\partial}G$, then $s_{1}- s_{2}$ 
would be a nonzero rational homomorphism from $G$ to its Lie algebra. 
By Condition (i) of Lemma 3.1, this gives the first part. If $(G,s)$ is a 
$D$-structure on $K$, so is $(G,\sigma(s))$ for any automorphism of the {\em field} ${\cal U}$ fixing $K$ pointwise.
So from uniqueness, $s$ is defined over $K$.
\newline
(ii) Let $\pi:\tilde B \to B$, and $\tau: T_{\partial}\tilde B \to \tilde B$. So 
$\pi\circ \tau:T_{\partial}\tilde B\to B$ has kernel a vector group, and thus there is a rational homomorphism
$s: \tilde B \to  T_{\partial}\tilde B$ such that  
\newline
(*) $\pi\circ\tau\circ s = \pi$. 
\newline
We claim that $s$ is a section of $\tau$, i.e.  that  $\tau\circ s = id$. Otherwise, by (*), $\tau\circ s - id$ is
a nontrivial rational homomorphism from $\tilde B$ to $ker(\pi)$. As the latter is a vector group, this contradicts $\tilde B$ being $asa$.
\newline
For part (iii) we simply quote Buium (\cite{Buium}, Theorem 3 of the Introduction). 
When $B$ is an abelian variety $A$ over 
$ \C(S)$,  this reflects the fact that by its very definition (cf. section A of the appendix), the class in $H^1(A, TA)$ of the  $TA$-torsor $T_\partial A$ is given by the Kodaira-Spencer map $\kappa(\partial)$, whose vanishing amounts to $A$ descending to $\C$.

\vspace{5mm}
\noindent
By virtue of Lemma 3.4 (i), we can and will talk about an almost semiabelian $D$-group $G$ without explicitly mentioning $s$.

\medskip
\begin{Lemma} Suppose $G$ is as $asa$ $D$-group. Then $G^{\partial}$ is the ``Kolchin closure" of $Tor(G)$, namely the smallest differential algebraic subgroup of $G$ containing $Tor(G)$.
\end{Lemma}
\pf  
Since the homomorphism $\partial \ell n_{G}$ takes values in a vector group, its kernel $G^{\partial}$ contains $Tor(G)$. Suppose for a contradiction that there is a proper differential algebraic subgroup $\cal H$ of $G^{\partial}$ containing $Tor(G)$.
By Fact 2.4 (ii), the Zariski closure $H$ of $\cal H$ would be a proper ($D$-)subgroup of $G$, contradicting the 
Zariski-denseness of $Tor(G)$ given by Lemma 3.1.iii.

\medskip
\begin{Corollary} If $G$ is an $asa$ $D$-group, then any rational homomorphism $f:G\to H$ from $G$ to a commutative algebraic $D$-group $(H,s)$ is a $D$-homomorphism. In particular any rational homomorphism between $asa$ $D$-groups is a $D$-homomorphism
\end{Corollary}
\pf We use Fact 2.4 (i). 
By the same argument as above,  the kernel $(H,s)^{\partial}$ of $\partial \ell n_{H,s}$ contains $Tor(H)$. Hence $f^{-1}((H,s)^{\partial})$, a differential algebraic subgroup of $G$, contains $Tor(G)$. By the
previous lemma, $f^{-1}((H,s)^{\partial})$ contains $G^{\partial}$, so $f(G^{\partial}) \subseteq H^{\partial}$. 

\vspace{5mm}
\noindent
Note the special case of Corollary 3.6 when $f$ is an embedding:
\begin{Corollary} If $H$ is an $asa$ $D$-group and is an algebraic subgroup of the commutative algebraic $D$-group $(G,s)$, then $H$ (with its unique $D$-group structure) is a $D$-subgroup of $(G,s)$.
\end{Corollary}

\vspace{5mm}
\noindent
Also Corollary 3.6 together with Lemma 3.1 and Lemma 3.4(ii)  yields:
\begin{Corollary} Let $G$ be an $asa$ $D$-group,
let $B$ be its maximal semiabelian quotient, and let $U$ be  a unipotent subgroup   of $\tilde B$  such that the  algebraic groups  $G$ and $ \tilde B/U$ are isomorphic. Then, $U$ is a $D$-subgroup of $\tilde B$. Note that if $G$ is defined over $K$, so are $\tilde B$ and $U$. 
\end{Corollary}

\vspace{5mm}
\noindent
Let us note in passing that the class of almost semiabelian $D$-groups is closed under quotienting by algebraic $D$-subgroups, but not of course under $D$-subgroups. Finally, we clearly obtain from Corollary 3.6   that if $B$ is a semiabelian variety, then  its universal vectorial extension $\tilde B$, equipped with its unique $D$-structure, is also universal in the category of  $D$-group extensions of $B$ by vector groups: namely if $(G,s)$ is
an algebraic $D$-group, and $G$ is, as an algebraic group,  an extension of $B$ by a vector group, then   then there is a unique morphism of $D$-groups from $ \tilde B$ to $ G$ satisfying the appropriate commutative diagram.

\vspace{8mm}
\noindent
{\em The $\sharp$-point functor on algebraic groups.}

\medskip

For the sake of completeness we tie up these notions with the {\em $\sharp$-point functor}, which will make a brief appearance in section 6 of the paper. 
Note that here, we do  not require $D$-group structures.
For $G$ an arbitrary connected commutative algebraic group  (over $\cal U$), 
we say that ``$G^{\sharp}$ exists" if  the intersection of all the Zariski-dense differential algebraic subgroups of $G$ is still Zariski dense in $G$, and we then denote it by $G^{\sharp}$.  Recall that if $G$ is almost semi-abelian (not necessarily a $D$-group) and if $B$ denotes  the semiabelian ``part" of $G$, then we have a canonical 
surjective homomorphism of algebraic groups $\pi:\tilde B \to G$, and moreover $\tilde B$ has a canonical $D$-group structure.

\begin{Proposition} Suppose that the algebraic group $G$ is almost semiabelian. 
\newline
(i) Then $G^{\sharp}$ exists and is the Kolchin closure of $Tor(G)$. It also equals $\pi({\tilde B}^{\partial})$.
\newline
(ii) Moreover, let $U_{1}$ be the maximal $D$-subgroup of $\tilde B$ contained in the (unipotent) subgroup $W := Ker(\pi)$.
Then $\pi$ induces an isomorphism of 
differential algebraic groups between $({\tilde B}/U_{1})^{\partial}$ and $G^{\sharp}$.
\end{Proposition}
\pf (i) 
By Lemma 4.2 of \cite{Pillay-countablemodels}, which is due originally to Buium, 
  any Zariski-dense differential algebraic subgroup of $G$ must contain $Tor(G)$, 
  but $Tor(G)$ is Zariski-dense. 
  As by 3.5, ${\tilde B}^{\partial}$ is the Kolchin closure of $Tor(\tilde B)$ and $\pi$ takes $Tor(\tilde B)$ onto $Tor(G)$, it is easy to conclude  the rest of part (i).
 \newline
  Notice that $G^\sharp = G^\partial$ when  we further assume that $G$ is a $D$-group; but in the general case under study here, the differential algebraic group $G^\sharp$ need not be defined by first order equations.
\newline
(ii) $W \cap {\tilde B}^{\partial}$ is a (connected) differential algebraic subgroup of ${\tilde B}^{\partial}$ so by
Fact 2.4 (ii), its Zariski closure is a (connected, unipotent) $D$-subgroup $U_{1}$ of $\tilde B$, and
$U_{1}^{\partial}$ is precisely  $W\cap {\tilde B}^{\partial}$. Then
(by the surjectivity of $\partial \ell n_{U_1}$) 
 $({\tilde B}/U_{1})^{\partial}$ is canonically isomorphic to ${\tilde B}^{\partial}/U_{1}^{\partial} = {\tilde B}^{\partial}/
W\cap {\tilde B}^{\partial}$ which is isomorphic to $G^{\sharp}$ under $\pi$.

\vspace{5mm}
\noindent
Note in particular that if $B$ is a semiabelian variety then $B^{\sharp}$ is canonically isomorphic to 
$(\tilde B/U_{B})^{\partial}$ where $U_{B}$ is the maximal unipotent $D$-subgroup of $\tilde B$.  
For instance, consider
the case where $B = A$ is a simple abelian variety defined over $K$, which does not descend to $C_K$ (equivalently, by simplicity, whose $C_K$-trace  is $0$), and let $U_A$ be the maximal unipotent $D$-subgroup  of $\tilde A$; so, $U_A$ is contained in  the kernel $W_{A}$ of $\tilde A \to A$. Since $A$ is not constant, Lemma 3.4.(iii) shows that $U_A$ is strictly contained in $W_A$. However, {\it $U_A$ need not vanish}, as witnessed by the  following example, which was shown to us by Y. Andr\'e: take a 
non constant
type IV abelian variety $A$ of even dimension $g = 2k \geq 4$, such that $\Q \otimes End(A)$ is a CM field of degree $g$, acting on $L\tilde A$ by  a CM type of the form $\{r_1 = s_1 = 1, r_2 = ...= r_k = 2, s_2 = ... = s_k = 0\}$. Then, the $g$-dimensional $K$-vector space $LW_A \simeq W_A$ is generated by two lines in $L\tilde A$ respectively contained in the planes where $F$ acts via the complex embedding  $ \sigma_1, \overline \sigma_1$, and by the planes where it acts via $ \sigma_2, ..., \sigma_k$, whereas $U_A$ is generated only by the latter planes, and therefore has dimension $g-2 > 0$ over $K$.

\vspace{5mm}
\noindent
\begin{Remark} If $B$ is a semiabelian variety over $K$, then it is convenient to have some notation for the algebraic $D$-group ${\tilde B}/U_{B}$. Let us call it $\overline B$.
With this notation, Proposition 3.9 (ii) gives a canonical isomorphism (over $K$) between 
${\overline B}^{\partial}$ and $B^{\sharp}$. The reader might think it natural to restrict our attention in Theorem 1.3 to $D$-groups of the the form ${\overline B}$. However,
our proofs will be of an inductive nature and involve taking $D$-quotients; just as with the discussion about $\tilde B$ after 
Remark 1.2,  this will of necessity force us into the general class of $asa$ $D$-groups.
\end{Remark}

\medskip
\subsection{Isoconstant $D$-groups}

One usually says that an algebraic group $G/K$ is isoconstant if there is an algebraic group $H/C_K$ such that $G$ and $H$ become isomorphic over an extension of $K$ . Since $K$ is algebraically closed, they then are automatically isomorphic over $K$.

\begin{Definition} Let $(G,s)$ be a commutative algebraic $D$-group defined over the algebraically closed differential field $K$. We will say that $(G,s_{G})$ is isoconstant if there is an algebraic group $(H,s_{H})$ such that $H$ is defined over ${\cal C}$ and $s_{H} = 0$,  and an isomorphism $f$ (over ${\cal U}$) of algebraic $D$-groups between $(G,s_{G})$ and $(H,s_{H})$.
\end{Definition}

So, the prefix  ``iso" here  refers to differential, rather than algebraic, extensions of the  (algebraically closed) base $K$. 
Once isoconstancy holds, we can actually insist that the data $(H,s_{H})$ and the isomorphism $f$ of this  definition be defined over $K^{diff}$, whereby  (as $C_{K} = C_{K^{diff}}$), $H$ will be over $C_{K}$. But again,  the isomorphism need not be defined over  $K$, as is shown by considering $D$-modules (cf. Remark 3.14 below). However,  for $asa$ $D$-groups, we do have rigidity:
\begin{Lemma} Suppose $G$ is an $asa$ $D$-group defined over the algebraically closed field  $K$, and suppose that $G$ is isoconstant. Then $G$ is isomorphic over $K$ to a constant $D$-group $(H,0)$.
\end{Lemma}
\pf By isoconstancy, and because $K$ is algebraically closed, the algebraic group $G$ is isomorphic over $K$ to an algebraic group $H$ defined over $C_{K}$. But then $H$ is $asa$, and by Corollary 3.6 this isomorphism is also one of $D$-groups.

\vspace{5mm}
\noindent
If $(G,s)$ is an algebraic $D$-group, it is not hard to see that it has a unique maximal connected isoconstant $D$-subgroup. 
For the next lemma we need to know that the image of any isoconstant algebraic $D$-group $(G,s)$ under a $D$-homomorphism $f$ is also isoconstant. This can be seen in various ways, one of which is as follows: We may assume $G$ defined over ${\cal C}$ and $s=0$. By Fact 2.3 (i), $Ker(f)$ is a $D$-subgroup of $G$, but $(Ker(f))^{\partial}$ is clearly a subgroup of $G({\cal C})$ and by Fact 2.4 (i) is Zariski-dense in $Ker(f)$. Hence $Ker(f)$ is defined over ${\cal C}$ and it is easy to conclude the argument.
\begin{Lemma} Let $G$ be an almost semi-abelian $D$-group defined over the algebraically closed field 
 $K$. Let $S$ be the maximal isoconstant connected $D$-subgroup of $G$, let $A$ be the abelian ``part" of $G$, and let  $U_{A}$ be the maximal unipotent $D$-subgroup of $\tilde A$. Then
\newline
(i)  $G/S$ has no toric part, and is therefore    a quotient of $\tilde A$.
\newline
(ii) If $S = \{0\}$, then  $A$ is an abelian variety with $C_{K}$-trace $0$, and  $G = \tilde A/U_{A}$.
\newline
(iii) Moreover all the objects  ($S$, $\tilde A$, etc.) and isomorphisms are defined over $K$.
\end{Lemma}
\pf  (i) The toric part of $G$ is by Corollary 3.7 a $D$-subgroup, so contained in $S$.
\newline
(ii) First of all $G$ has no toric part, so is a quotient of $\tilde A$ by a unipotent $D$-subgroup $U$. As any commutative unipotent $D$-group is isoconstant, our assumption  on $G$ forces $U$ to  be the maximal unipotent $D$-subgroup of $\tilde A$. Now, we may work up to isogeny, and  assume $A$ to be of the form $A_{0}\times A_{1}$ where $A_{0}$ is over $C_{K}$ and $A_{1}$ has $C_{K}$-trace $0$. Then ${\tilde A} =
{\tilde A_{0}} \times {\tilde A_{1}}$. But $\tilde A_{0}$ is a constant $D$-group, as is its image 
 in $G$ (by the paragraph preceding this lemma). So $A_{0} = 0$, and we see that $G$ must be the quotient of $\tilde A_{1}$ by its maximal unipotent $D$-subgroup, as required.
\newline
(iii) $S$ is defined over $K$ by uniqueness. The rest is clear.

\medskip
\noindent
\begin{Remark}  Let $(G,s)$ be a (commutative) unipotent $D$-group over $K$. Then $(G,s)$ is isoconstant.
\end{Remark}
\pf  This of course belongs to the theory of linear differential equations. Identifying $G$ with $\G_{a}^{n}$, $s$ has the form $x\mapsto (x,Ax)$ for $A$ an $n\times n$ matrix over $K$. We can find a basis $\{v_{1},..,v_{n}\}$ of $G$ (as a vector space over $\cal U$) which is simultaneously a basis of the ${\cal C}$-vector space $(G,s)^{\partial}$. Now,
$v_{1},..,v_{n}$ can be chosen from $K^{diff}$ (but not always from $K$) and generate the Picard-Vessiot extension 
for the linear differential system $\partial(-) = A(-)$. With respect to the basis $\{v_{1},..,v_{n}\}$, $s$ becomes $0$, and so over $K(v_{1},..,v_{n})$, $(G,s)$ is isomorphic to $(\G_{a}^{n},0)$.

\section{Further ingredients and special cases}

Here we present the key ingredients for the proof of Theorem 1.3. 

\subsection{The socle theorem}
Proposition 4.1 below follows from Corollary 3.11 of \cite{Pillay-Ziegler} (also appearing in the language of $D$-groups as Corollary 2.11 in 
\cite{Kowalski-Pillay}). The latter generalizes and has its origin in what is often called the ``socle theorem" of Hrushovski. This ``socle theorem" is actually a combination of Proposition 4.3 \cite{Hrushovski} with the validity of the Zilber dichotomy in the theory of differentially closed fields of characteristic $0$.  In fact in the case at 
hand, what is needed can  probably  be extracted from Hrushovski's results. In any case the relevant statement concerns 
commutative connected finite-dimensional differential 
 algebraic groups, equivalently commutative, connected, 
groups of finite Morley rank ${\cal B}$ 
definable in the differentially closed field ${\cal U}$. By the {\em algebraic 
socle}, $as({\cal B})$, of ${\cal B}$ we mean the maximal connected definable subgroup of ${\cal B}$ which is definably isomorphic to a 
group of the form $C({\cal C})$ for $C$ some algebraic group over ${\cal C}$. If $X$ is a differential algebraic 
subset of ${\cal B}$, its stabilizer $Stab_{{\cal B}}(X)$ is defined to be $\{g\in {\cal B}: g+ X = X\}$. The result says that 

\vspace{2mm}
\noindent
(*) if $X$ is an 
irreducible differential algebraic subset of ${\cal B}$ such that $Stab_{{\cal B}}(X)$ is finite, then some translate of $X$ is 
contained in $as({\cal B})$. 

\vspace{2mm}
\noindent

It should be mentioned that this result is quite powerful, and together with the kind of material in section 3, yields quickly a proof of the function field Mordell-Lang conjecture in characteristic $0$.
Now if ${\cal B} = G^{\partial}$ for $G$ an algebraic $D$-group, then
it follows from Fact 2.4(ii) that $as({\cal B}) = H^{\partial}$ where $H$ is the maximal connected isoconstant $D$-subgroup of $G$,
as introduced before Lemma 3.13. Bearing this in mind, Proposition 4.1 below is simply a principal homogeneous space version of 
the socle theorem
(*) above.

\begin{Proposition} Let $G$ be an almost semiabelian $D$ group over $K$.
Let $Y\subset G$ be a translate (coset) of $G^{\partial}$ and $Z$ an irreducible differential
algebraic subset of $Y$. Let 
${\cal S}< G^{\partial}$ be the stabilizer of $Z$ (with respect to the action
of $G^{\partial}$ on $Y$). Suppose that ${\cal S}$ is finite. Then $Z$ is contained in a coset $Z'$ of 
$H^{\partial}$ where $H$ is the maximal isoconstant $D$-subgroup of $G$. Moreover
if $Y,Z$ are defined over $K$, so is $Z'$.
\end{Proposition}

\subsection{On the semi-constant part and descent of $({\bf HX})_{K}$}
$K$ will  
again denote an algebraically closed differential field.
Suppose 
 $B$ is a semiabelian variety over $K$, with
$0 \mapright{} T \mapright{}
     B \mapright{\pi} A \mapright{}
      0 ~$  the canonical short exact sequence.
Let $A_{0}<A$ be the $C_{K}$-trace of $A$  (maximal abelian subvariety of $A$ isomorphic to an abelian variety over $C_{K}$),
 and $B_{0} = \pi^{-1}(A_{0})$. 
 As in the introduction, 
 we call $B_{0}$ the semi-constant part of 
 $B$. More generally,
suppose $G$ is an almost semi-abelian variety with semiabelian part $B$; in particular, $G = {\tilde B}/U$ for   some unipotent subgroup $U$ of $\tilde B$. We denote $B$ by $G^{sa}$ and define the {\it semi-constant part of $G$}  to be $G^{sa}_{0}$, namely $B_{0}$. 
In this setting, our hypothesis  $({\bf HG})_{0}$ on the almost semi-abelian group $G$ can be stated as follows:

\vspace{2mm}
\noindent
$({\bf HG})_{0}$  {\it the semi-constant part $B_{0} = G^{sa}_{0}$ of $G$ is an isoconstant algebraic group.}

\vspace{2mm}
\noindent 
In other words, the semi-abelian variety $G^{sa}_{0}$ is  isomorphic (and then, automatically over $K$) to an algebraic group defined over $C_{K}$.
It is clear that this hypothesis $({\bf HG})_{0}$ is 
preserved under quotients,

\medskip
We will describe concretely $(LG)^{\partial}(K)$ under the assumption $({\bf HG})_{0}$ and conclude that the hypothesis $({\bf HX})_{K}$ is 
then
preserved under quotients.

\medskip
From now on, we assume that the algebraically closed differential field $K$ has transcendence degree $1$ over its field $\C$ of constants. We will make use of two classical results from Deligne's  Hodge II \cite{Deligne} which we now explain. Let $A$ be an abelian variety over $K$, and $\tilde A$ its universal vectorial extension, equipped with its unique $D$-group structure. From subsection 2.3 this provides $L\tilde A$ with a connection $\partial_{L\tilde A}$. In section H of the appendix,  we check that this connection is the dual of the Gauss-Manin connection on $H^{1}_{dR}(A/K)$. 
Write $A$ up to isogeny as $A_{0}\times A_{1}$ where $A_{0}$ descends to $\C$ and $A_{1}$ as $\C$-trace $0$. Then $\tilde A = {\tilde A_{0}} \times {\tilde A_{1}}$. We will assume that $A_{0}$ is already over $\C$, thus so is ${\tilde A_{0}}$. As the $D$-group structure on ${\tilde A_{0}}$ is trivial, 
 so is the corresponding connection on $L{\tilde A_{0}}$. Consequently, $(L{\tilde A_{0}})^{\partial} = L{\tilde A_{0}}(\C)$. The above mentioned results  
translate  in the present setting as follows (see   Corollary H.5 of the appendix):
\newline
(I) (Semisimplicity) The $D$-module  $L{\tilde A}$ is semisimple, i.e. a direct sum of simple $D$-modules over the algebraically closed field $K$ (Deligne, cf. \cite{Deligne}, II, 4.2.6). 
\newline
(II) (Theorem of the fixed part)  An horizontal vector of ${\partial_{L\tilde A}}$ which is invariant under the monodromy lies in $L\tilde A_0$ (Griffiths-Schmid-Deligne, cf. \cite{Deligne}, II, 4.1.2). In other words,   $(L{\tilde A})^{\partial}(K) = L{\tilde A_{0}}(\C)$.

\begin{Lemma} (i) Let $B$ be a semiabelian variety over $K$. Assume that its semi-constant part $B_{0}$ is defined over $\C$, and let $\tilde B_{0} \subseteq \tilde B$ be the universal vectorial extension of $B_{0}$. Let $U$ be a unipotent $D$-subgroup of $\tilde B$,  let $G$ be the $asa$ $D$-group $\tilde B/U$, and let $U_0 =   \tilde B_{0} \cap U$. Then,
$(LG)^{\partial}(K) = L(\tilde B_{0}/U_0)(\C)$.
\newline
(ii) Let $G$ be an $asa$ $D$-group over $K$ satisfying $({\bf HG})_{0}$. Let $x\in LG(K)$ be such that for no proper algebraic subgroup $H$ of $G$ over $K$ is $x\in LH(K) + (LG)^{\partial}(K)$. Let $G'$ be
a quotient of $G$ by a $D$-subgroup defined over $K$, so itself an $asa$ $D$-group. Let $x'\in LG'(K)$ be the image of $x$ under the corresponding projection $LG \to LG'$. Then for no proper algebraic subgroup $H'$ of $G'$ over $K$ is $x'\in LH'(K) + (LG')^{\partial}(K)$.
\end{Lemma} 

\pf  (i) Note first that $U_0 $ is a $D$-subgroup of the constant $D$-group $\tilde B_{0}$ so is also defined over the constants and has trivial $D$-group structure. Hence $(L(\tilde B_{0}/U_0))^{\partial}$ is the set of constant points of $L(\tilde B_{0}/U_0)$.  Let $T$ and $A = A_{0}\times A_{1}$ be as usual for $B$. Identifying $U$ with a unipotent $D$-subgroup of $\tilde A$, we have the exact sequence
\newline
$0 \to LT(\C) \to (LG)^{\partial}(K^{diff}) \to (L(\tilde A/U))^{\partial}(K^{diff}) \to 0$.
\newline
Let $x\in (LG)^{\partial}(K)$, and let $\overline x$ be its image in $(L(\tilde A/U))^{\partial}(K)$. By (I) above, $\overline x$ lifts to a point  $\overline \xi \in (L\tilde A)^{\partial}(K)$. By (II) above the latter space is precisely
$(L\tilde A_{0})(\C)$. Hence ${\overline x}\in L(\tilde A_{0}/U_0)(\C)$. This now implies that $\overline x$ lifts to a point $\xi \in L(\tilde B_{0}/U_0)(\C)$. But $x,\xi \in (LG)^{\partial}(K)$ and have the same image in 
$(L(\tilde A/U))^{\partial}$, hence their difference lies in $LT(\C)$, which is contained in 
$L(\tilde B_{0}/U_0)(\C)$. Thus,  already $x\in L(\tilde B_{0}/U_0)(\C)$, and (i) is proved.

\medskip
\noindent
(ii) Let $G = \tilde B/U$ and  $G' = \tilde B'/U'$  (with $B, B'$ semiabelian and $U,  U'$ unipotent $D$-subgroups of $\tilde B$, $\tilde B'$ respectively). The hypothesis $({\bf HG})_{0}$ says that the semi-constant part $B_0$ of $B$  is isoconstant. As $B'$ is a quotient of $B$, the semiconstant part $B'_0$ of $B'$ is also isoconstant. 
Assume, as we may,  that $B_{0}$ is already defined over $\C$, and set  $U_0 = \tilde B_0 \cap U, U'_0 = \tilde B'_0 \cap U'$.
It is clear that under the relevant quotient map $f:G \to G'$, $\tilde B_{0}/U_0$ maps onto
$\tilde B'_{0}/U'_0$. 
As both the latter groups are defined over $\C$, we see that
$f$ maps $(\tilde B_{0}/U_0)(\C)$ onto $(\tilde B'_{0}/U'_0)(\C)$. By part (i) of this lemma, it follows that $f$ maps
$(LG)^{\partial}(K)$ onto $(LG')^{\partial}(K)$. 

\vspace{2mm}
\noindent
Now let $x,x'$ be as in the statement of part (ii) of the lemma, and suppose towards a contradiction, that
$x' = u' + z'$ with $u' \in LH'(K)$ for some proper algebraic subgroup $H'$ of $G'$ defined over $K$   and $z' \in (LG')^{\partial}(K)$. By what we have just observed, $z'$ lifts to a point $z\in (LG)^{\partial}(K)$.
Then $x-z \in LG(K)$ and $f(x-z) = u'\in LH'(K)$. Hence $x-z \in LH(K)$ where $H := f^{-1}(H')$ is a proper algebraic subgroup of $G$ defined over $K$, contradicting our hypothesis on $x$. Lemma 4.2 is proved.

\subsection{Special cases: Ax and the theorem of the kernel}
We give two special cases where Theorem 1.3 holds (and where the $({\bf HG})_{0}$ is not mentioned but it automatically holds). The first is a special case of Ax's theorem, as slightly generalized in \cite{Bertrand}, although there are other direct proofs 
(see, e.g. \cite{Kirby}). Again $K$ is an algebraically closed differential field of transcendence degree $1$ over its field $\C$ of constants.
\begin{Proposition} Suppose that $G$ is an isocontant $asa$ $D$-group over $K$, $x\in LG(K)$, and $y\in G({\cal U})$ such that $\partial\ell n_{G}(y) = \partial_{LG}(x)$. Assume that $x\notin LH(K) + (LG)^{\partial}(K)$ for any proper algebraic subgroup of $G$ defined over $K$. Then $tr.deg(K(y)/K) = dim(G)$.
\end{Proposition}
\pf By Lemma 3.12 we may assume that $G$ is defined over $\C$, in which case $\partial\ell n_G$ is Kolchin's logarithmic derivative, $\partial_{LG}$ is just $\partial: LG \to LG$ and so $(LG)^{\partial}(K^{diff}) = LG(\C)$.
The semiabelian part $B$ of $G$ is defined over $\C$, and  $G = \tilde B/U$ for   some unipotent subgroup $U$ of $\tilde B$ over $\C$

We may assume that $y$ lies in  $G(K^{diff})$.   Let $\tilde x \in L{\tilde B}(K)$ be a lift of $x$, let $\tilde y\in \tilde B(K^{diff})$ be a solution of
 $\partial\ell n_{\tilde B}(-) = 
 \partial_{L\tilde B} \tilde x$ 
and let $\overline y$ be the projection of $y$   to $B$. We first check that the relative hull $B_{\overline y}$ of $\overline y$ fills up  $B$, i.e. that there is no proper algebraic subgroup $\overline H$ of $B$ (over $\C$) such that $\overline y \in \overline H({\cal U}) + B(\C)$. Otherwise, the projection $\overline x$ of $x$ to $LB(K)$ would satisfy
$\partial_{LB}(\overline x ) = \partial \ell n_B (\overline y) \in L\overline H({\cal U}) $, so
 $\overline x \in \big(L\overline H({\cal U}) + LB(\C)\big) \cap LB(K) = L\overline H(K) + LB(\C)$, and if $H$ denotes the inverse image of $\overline H$ in $G$, $x$ would lie in $LH(K) + LG(\C)$, contradicting  $({\bf HX})_{K}$.

Now by Proposition 1b of \cite{Bertrand},  $tr.deg(K(\tilde y)/K) = dim({\tilde B})$. Let $y'$ be the image of ${\tilde y}$ under the canonical projection from ${\tilde B}$ to $G$. Then, $tr.deg(K(y')/K) = dim(G)$. 
Now both $y,y'$ are in $G(K^{diff})$, and $\partial\ell n_{G} (y) = \partial \ell n_{G}(y') = \partial_{LG}x$. Hence 
$y-y' \in Ker(\partial\ell n_{G})(K^{diff}) = G(\C)$. So $tr.deg(K(y)/K) = dim(G)$ as required.

\vspace{5mm}
\noindent
The next special case is fundamental and depends on Manin's theorem of the kernel, in a stronger form due to Chai \cite{Chai}. 
 A direct proof of this result is given in  section K of the Appendix 
 (see Remark K.2 for a discussion of Chai's full sharpening). We take $G$ to be an almost abelian $D$-group over $K$, namely for some abelian variety $A$ over $K$, $G$ is a quotient of $\tilde A$ by  
 some unipotent $D$-subgroup $V$. As usual, we have $0\to W_{A} \to \tilde A \to A \to 0$ over $K$, and we denote by $U_A$ the  maximal unipotent $D$-subgroup of $W_{A}$. So,  $G = {\tilde A}/V$, with $V \subset U_A$.    
 It follows that $W_{A}/V$ is the unipotent part of $G$ and we write it as $W_{G}$. 
$W_{G}$ can be identified with its Lie algebra $LW_{G}$, which is itself contained in $LG$.

With the above notations, Chai's theorem reads as follows (cf. Theorem K.1):
\begin{Proposition}  Suppose $A$ has $\C$-trace $0$. Suppose that $x\in LG(K)$ and $y\in G(K)$ are such that
$\partial\ell n_{G}(y) = \partial_{LG}(x)$. Then $x\in LW_{G}(K)$.
\end{Proposition}

Notice that conversely, given $x \in LW_G(K)$, there does exist a point $y \in G(K)$ such that $\partial\ell n_{G}(y) = \partial_{LG}(x)$, namely $ y = x$ itself, viewed as a point in $W_G \subset G$. Indeed, Corollary G.4 of the Appendix shows that $\partial\ell n_{G}$ and $\partial_{LG} $ coincide on $W_G \simeq LW_G$.

\section{Proofs of main results}
Throughout this section $K$ is an  algebraically closed differential field of transcendence degree $1$ over its field of constants which is assumed to be $\C$.

\subsection{Proof of Theorem 1.3}

We recall that $G$ is an almost semi-abelian $D$-group over $K$ satisfying $({\bf HG})_{0}$  (its semiconstant part is constant). We take $x\in LG(K)$ satisfying $({\bf HX})_{K}$: $x\notin L(H)(K) + (LG)^{\partial}(K)$ for any proper algebraic subgroup $H$ of $G$ over $K$.  Given any $y\in G({\cal U})$ satisfying $\partial\ell n_{G}(y) = \partial_{LG}(x)$, we must prove that $tr.deg(K(y)/K) = dim(G)$.

\bigskip
Our proof is of a differential algebraic nature. We work in ${\cal U}$, that is we identify (differential) algebraic varieties and groups with their pointsets in ${\cal U}$. We first recall the notion of ``generic points" of differential algebraic varieties.  
For $X$ an irreducible finite-dimensional differential algebraic variety, defined over $K$, by a generic point of $X$ over $K$ we mean a point $\alpha \in X$ 
 such that $tr.deg(K\langle \alpha \rangle/K)$ is maximum possible. If $X$ is of the form $(G,s)^{\partial}$ for some algebraic $D$-group $(G,s)$ defined over $K$, then for $\alpha \in G^{\partial}$, $K\langle \alpha \rangle = K(\alpha)$ and clearly the maximum possible $tr.deg(K(\alpha)/K)$ equals $dim(G)$ (as $G^{\partial}$ is   a differential algebraic subvariety of $G$ which is Zariski-dense in $G$). 
 Likewise if $X$ is a translate of $G^{\partial}$ inside $G$ and is defined over $K$, then
max.$\{tr.deg(K\langle \alpha \rangle/K): \alpha\in X\} = dim(G)$.

\medskip
\noindent
We now begin the proof proper.
\newline
Let $dim(G) = n$. Let $a = \partial_{LG}(x) \in LG(K)$, and let $Y\subseteq G$ be the solution set of
$\partial\ell n_{G}(-) = a$. So $Y$ is a coset (translate) of $G^{\partial} = Ker(\partial\ell n_{G})$ in $G$, defined over $K$. By the discussion of genericity above, for generic $y\in Y$ (over $K$), $tr.deg(K(y)/K) = n$. 
 Our desired conclusion is that {\em for all} $y\in Y$, $tr.deg(K(y)/K) = n$. We claim that this is 
{\em equivalent} to saying that $Y$ has no proper differential algebraic subsets defined over $K$. 
Let us explain briefly the equivalence. The differential equations defining $Y$ in $G$ give, for all $y \in Y$, a rational control of $\partial y$ in terms of $y$.  
Hence differential algebraic subvarieties of $Y$ defined over $K$ have the form $Z\cap Y$ for $Z$ algebraic subvarieties of $G$ defined over $K$, and our desired conclusion that $Y$ meets no proper algebraic subvariety of $G$ defined over $K$ is equivalent to $Y$ having no differential algebraic subvariety defined over $K$. (This is is course related to a version of Fact 2.4 (ii) for subvarieties rather than algebraic subgroups.)

So we will prove:

\begin{Lemma} $Y$ has no proper irreducible differential algebraic subset defined over $K$. (Or
in model-theoretic language the formula ``$y\in Y$" isolates a complete type over $K$.)
\end{Lemma} 

\pf  Globally the proof proceeds by induction on $n$. If $n=1$ then $G$ is either $\G_{m}$ or an elliptic curve $E$.
In either case (owing to Lemma 3.4 (iii)), $G$ descends to $\C$ and so we finish by Proposition 4.3.

So let us assume  $n>1$. We will suppose that $Z$ is a proper irreducible differential algebraic subset of $Y$, defined over $K$ and look for a contradiction.
This will be a somewhat involved case analysis, reducing to the special cases discussed in section 4.3.

Let ${\cal S} := \{g \in G^{\partial}: g+Z = Z\} < G^{\partial}$ be the stabilizer of $Z$. 
 Then $\cal S$ is a differential algebraic subgroup of $G^{\partial}$ defined over 
$K$. By Fact 2.4 (ii), ${\cal S} = S\cap G^{\partial}$, where $S$ is the Zariski closure of $\cal S$ 
 and is an algebraic $D$-subgroup of $G$ defined over $K$. Note that $S$ is a proper subgroup of $G$, for otherwise ${\cal S} = G^{\partial}$ and $Z = Y$. 

\medskip
\noindent
{\bf CASE I.}  $\cal S$ is infinite.
\newline
Let $G'$ 
 be $G/S$ (in fact it will be enough to quotient by the connected component of $S$). 
Let $\pi:G\to G'$ be the canonical $K$-rational surjective homomorphism. So $G'$ is an
almost semiabelian $D$-group, defined over $K$, with dimension positive and $< n$ and $\pi$ is a homomorphism
of $D$-groups, inducing a surjective $K$-rational homomorphism  $L\pi$ from $LG \to LG'$, 
which is also a homomorphism of $D$-modules by 2.5. It follows that $\pi(Y)$ is the solution space of the equation 
$\partial\ell n_{G'}(-) = a'$ where $a' = L\pi(a)$. Note also that  $\partial_{LG'}(x') = a'$ where $x' = 
L\pi(x) \in LG'(K)$. By Lemma 4.2 the hypothesis $({\bf HX})_{K}$ holds for $x'$. By the induction hypothesis, $\pi(Y)$ has no proper differential algebraic subset defined over $K$. As $\pi(Z)$ is a differential constructible subset of $\pi(Y)$, defined over $K$, it follows that $\pi(Z) = \pi(Y)$. This implies that $Y\subseteq Z + S$. As $Z\subseteq Y$, $Y$ is a $PHS$ for $G^{\partial}$ and ${\cal S} = S\cap G^{\partial}$,  it follows that $Y = Z + {\cal S}$. This contradicts ${\cal S}$ being the stabilizer of $Z$ and $Z$ being a proper differential algebraic subset of $Y$. So CASE I leads to a contradiction.

\medskip
\noindent
{\bf CASE II.} $\cal S$ is finite. 
\newline
Let $H$  denote the maximal connected isoconstant $D$-subgroup of $G$ (which is defined over $K$). By
Proposition 4.1, $Z$ is contained in a coset (i.e. orbit) $Z'$ of $H^{\partial}$, and $Z'$ is defined over $K$.
We have again a $K$-rational surjective homomorphism $\pi:G \to G/H$ of $D$-groups, with $L\pi:LG\to L(G/H)$ a surjective homomorphism of $D$-modules. And again, $L\pi(x) = x' \in L(G/H)(K)$, and 
$\partial_{L(G/H)}(x') = a' = L\pi(a) \in L(G/H)(K)$. But now 
 $\pi(Z) = y'$ is a {\it point} in $\pi(Y) = \{$ solutions of $\partial\ell n_{G/H}(-) = a' \}$ since $Z \subset H$, and this point is $K$-{\it rational} since  $\pi(Z) = \pi(Z')$ and $Z'$ is defined over $K$.
Moreover the hypotheses $({\bf HG})_{0}$ and $({\bf HX})_{K}$ remain valid for $G/H$ and $x'$ (using Lemma 4.2).

\medskip
\noindent
We have three subcases:
\newline
(a) $H = G$. 
\newline
This means that $G$ itself is an isoconstant $D$-group, and we contradict Proposition 4.3.

\vspace{2mm}
\noindent
(b) $H$ is a proper non-zero 
subgroup of $G$.
\newline  So $dim(G/H) $ is both positive and $< n$,  and we can use the induction hypothesis. Since  
$tr.deg(K(y')/K) = 0$, we 
have a contradiction.

\vspace{2mm}
\noindent
(c) $H = \{0\}$ . 
\newline
By Lemma 3.13 (ii) $G = {\tilde A}/U_{A}$, where $A$ is an abelian variety over $K$ with $\C$-trace $0$ and $U_{A}$ is the maximal unipotent $D$-subgroup of ${\tilde A}$. Moreover $y = y'\in G(K)$ and $x\in LG(K)$ satisfy
$\partial\ell n _{G}(y) = \partial_{LG}(x)$. 
By Proposition 4.4, $x\in W_{G}(K) = LW_{G}(K)$, where $W_{G}$ is a proper algebraic subgroup of $G$ defined over $K$. This contradicts the 
hypothesis $({\bf HX})_{K}$.

\bigskip
\noindent
We have shown that all cases lead to a contradiction. So Lemma 5.1 is proved, as is Theorem 1.3.

\subsection{Dropping $({\bf HG})_{0}$}
We give a promised version of Theorem 1.3.
\begin{Theorem} Let $G$ be an $asa$ $D$-group over $K$, $x\in LG(K)$ and assume that
\newline
$({\bf HX})$: for no proper algebraic subgroup $H$ of $G$ defined over $K$ is 
$x\in LH(K^{diff}) + (LG)^{\partial}(K^{diff})$.
\newline
Then for any $y\in G 
({\cal U})$ 
such that $\partial\ell n_{G}(y) = \partial_{LG}(x)$, $tr.deg(K(y)/y) = dim(G)$.
\end{Theorem}

\pf  So note that there is now no restriction on the semiconstant part of $G$. But the 
``arithmetic" hypothesis $({\bf HX})_{K}$ on $x$ has been strengthened to the differential algebraic
hypothesis ${\bf HX}$.  The proof is identical to that of Theorem 1.3 above, except that 
the stronger hypothesis $({\bf HX})$ is easily seen to descend under quotients by $D$-subgroups, while as already noted, there is no $({\bf HG})_{0}$ assumption in the special cases covered by Lemmas 4.3 and 4.4 (in fact it is automatically true there.)

\medskip

We take opportunity of this discussion to point out that in  
either condition $({\bf HX})_{K}$ or  $({\bf HX})$, 
we must consider all proper {\it algebraic} subgroups $H$ of $G$, not just its algebraic $D$-subgroups. For instance,  if $G = \tilde A$ for some non constant elliptic curve  $A/K$, so that $(LG)^\partial(K) = 0$,  and if $0 \neq x \in LW_A(K)$, then, $y = x$, viewed as a point of $W_A \subset G$, satisfies $\partial \ell n_G y = \partial_{LG} x $, in view of Corollary G.4 of the Appendix (see also the remark following Proposition 4.4). We then have  $tr.deg. (K(y)/K) = 0 < dim(G)$, although   
since $U_A$ here vanishes,  $x$ lies in the Lie algebra of no  proper $D$-subgroup of $G$.

\subsection{A 
counterexample}
We give the simplest possible example showing that in Theorem 1.3, the $({\bf HG})_{0}$ hypothesis cannot be dropped in general. 
Let 
$K= \C(z)^{alg}$, let $  E$ be an elliptic curve defined  over $\C$ and let $B$  be a nonconstant extension of $E$ by $\G_{m}$, defined over $K$ (such extensions are given by $K$-rational points on the dual $\hat E$ of $E$, not lying in $\hat E(\C)$). We take as our $asa$ $D$-group the universal vectorial extension $G = \tilde B$ of $B$, and recall from Corollary 3.7 (or from Fact H.3 of the Appendix) that   $LG$ is an
 extension of $L\tilde E$ by $L\G_m$ in the category of  $D$-modules over $K$.  

Let ${\overline x}$ be a nonzero point in $LE(\C)$, which we  lift to a point $\tilde x \in L{\tilde E}(\C)$  
and  finally to a $K$-rational point 
$x\in LG(K)$ of $LG$. Then:

\medskip
\noindent
{\em Claim I.}  $  G$ does not satisfy $({\bf HG})_{0}$. Indeed, its semi-constant part $G^{sa}_{0}$ is  
$B$ itself, which is not isoconstant.

\medskip
\noindent
{\em Claim II.} If $y\in { G}(K^{diff})$ satisfies $\partial\ell n_{  G}(y) = \partial_{L  G}(x)$, then $tr.deg(K(y)/K) \leq 1$.
\newline
\pf  Let $\partial_{L G}(x) = a \in LG(K)$. Then the image of $a$ under the projection to
$L{\tilde E}$ is  $\partial_{L\tilde E}({\tilde x}) = 0$. Hence the solution set of 
$\partial\ell n_{  G}(-) = a$ in ${  G}(K^{diff})$ projects onto the solution set
of $\partial\ell n_{\tilde E}(-) = 0$ in ${\tilde E}(K^{diff})$ which is precisely ${\tilde E}(\C)$.  Since the fibers of the projection $G \to \tilde E$ are one-dimensional, the claim follows.

\medskip
\noindent
{\em Claim III.}  $(L G)^{\partial}(K) =  
 (L\G_{m})^{\partial}(K)$  (= $L\G_{m}(\C)$). 
\newline
\pf   Since $\G_m$ is the maximal constant subgroup of $B$,  
this follows from the extension to mixed Hodge structures of the theorem of the fixed part, as given in \cite{Steenbrink-Zucker}, Prop.  4.19.

\medskip
\noindent
{\em Claim IV.}  $x$ satisfies $({\bf HX})_{K}$.
\newline
\pf  The only proper algebraic subgroups $H$  of $ G$ are $\G_{m}$, $\G_{a}$ and $\G_{m}\times \G_{a}$.
So if $x\in LH(K) + (LG)^{\partial}(K)$ then by Claim III, $x\in L(\G_{m}\times \G_{a})$ so could not project to a nonzero element of $LE$. This contradicts the choice of $\overline x$. 

\medskip
\noindent
Claims II and IV show that Theorem 1.3 is in general false without the $({\bf HG})_{0}$ hypothesis. 
On the other hand, Theorem 5.2 remains valid: since $\tilde x \in (L\tilde E)^\partial$, the point $x$ in this example lies in $L\G_m({K^{diff}}) + (LG)^\partial({K^{diff}})$, and therefore violates Hypothesis $({\bf HX})$.

\section{$K$-largeness and a differential Galois-theoretic proof}
We will make an additional assumption, $K$-largeness, on our $asa$ $D$-group $G$ over $K$, and obtain the conclusion
of Theorem 1.3, replacing the use of the ``socle theorem" (Proposition 4.1)  by rather softer tools. At the time of writing we are aware that the $K$-largeness hypothesis is very strong, but nevertheless as such Galois-theoretic methods were 
the initial motivation of this paper (cf. \cite{Bertrand}, Remark 1), it seems appropriate to include this material.  
More to the point, the proof of Theorem 1.3 which these methods provide in this case is the exact copy of Kolchin's classical proof of the Ostrowski theorem and of its multiplicative analogue, cf. Section 2 of  \cite{Kolchin-Ostrowski}, and this suffices to justify their insertion here.
 
\bigskip
\subsection{$K$-large algebraic $D$-groups}

We begin in our general context, where $K$ is an arbitrary  (small) algebraically closed differential subfield of ${\cal U}$, all in characteristic 0. For $(G,s)$  an algebraic $D$-group defined over $K$, we say, following \cite{Pillay}, that {\it $(G,s)$ is $K$-large if 
$G^{\partial}(K^{diff}) = G^{\partial}(K)$}. This makes sense for arbitrary $D$-groups, but we will restrict to the case where $G$ is commutative and use additive notation.

\medskip

Here are some examples, still with $K$  arbitrary. A constant algebraic $D$-group $(G, s = 0)$ is automatically $K$-large (and even $C_K$-large). 
A unipotent $D$-group, although isoconstant, is $K$-large if and only if the corresponding  $D$-module  is completely solvable over $K$, while if $A$  
is an abelian variety over $K$ 
 with $\overline A = \tilde A/U_A$ as in Remark 3.10, 
  $T$ is  a torus, and  $G$ is of the form $T\times{\overline A}$, then $G$ is $K$-large. In other words, the sharp points of an abelian variety $A$ always satisfy $A^\sharp(K^{diff}) = A^\sharp(K)$.  
 For the latter fact, it is enough to consider a simple abelian variety $A$, in which case the statement that 
$A^{\sharp}(K^{diff}) = A^{\sharp}(K)$ is precisely Lemma 2.2 of \cite{Marker-Pillay} (which depends on work of Hrushovski and Sokolovic). 
On the other hand, it follows from Proposition 6.1 below that over $K = \C(z)^{alg}$, the algebraic $D$-group $G/K$ considered in subsection 5.3 is not $K$-large.

\bigskip
We now establish a relationship between the hypothesis $({\bf HG})_{0}$ and $K$-largeness, in a quite general setting, using model-theoretic methods. Recall from Remark 3.10 
 the notation $\overline B$ for $B$ a semiabelian variety:  $\overline B$ is the quotient of $\tilde B$ by its maximal unipotent $D$-subgroup $U_B$.

\begin{Proposition} Let $K$ be an (algebraically closed) differential field. Assume that $C_{K}$ has infinite transcendence degree and that $K$ is the algebraic closure of a differential field which is finitely generated (as a differential field) over $C_{K}$. Let $B$ be a {\rm  semi-constant} 
 semiabelian variety over $K$,  
 i.e. such that the abelian part of $B$ descends to $C_{K}$. Then the following are
equivalent:
\newline
(i) $B$ is constant, i.e. descends to $C_{K}$.
\newline
(ii) The algebraic $D$-group $G = \overline{B}$ is $K$-large.
\end{Proposition}
\pf (i) implies (ii) is immediate: We may assume that $B$ is defined over $C_{K}$. But then $\tilde B$ is defined over $C_{K}$ and its unique $D$-group structure is the trivial one. Hence $U_{B} $ coincides with the maximal unipotent subgroup $W_{B}$ of $\tilde B$, and $G := \tilde B/U_B = B$ with the trivial $D$-group structure.   
But then $G^{\partial}(K^{diff}) = G^\partial(C_{K^{diff}}) = G(C_{K}) = G^{\partial}(K)$.
\newline
(ii) implies (i). There is no harm in assuming that both $T$ and $A$ are defined over $C_{K}$ (where $T, A$ are the toric, respectively abelian, parts of $B$). 
We have our exact sequence $0 \to T\to {\tilde B} \to {\tilde A} \to 0$ of $D$-groups, and note that ${\tilde A}$ is defined
over $C_{K}$ and has trivial $D$-structure, while
$G = \overline B$ is obtained by quotienting $\tilde B$ by its maximal unipotent $D$-subgroup  $U_B$. Note that $U_B$ identifies with a unipotent $D$-subgroup of $\tilde A$ 
 and is thus defined over $C_{K}$.
 Then $A^{\flat} := {\tilde A}/U_B$  is also defined over $C_{K}$, has trivial $D$-structure, and sits in the exact sequence of $D$-groups
\newline 
(\dag)  $0 \to T\to   {G} \to A^{\flat} \to 0$.

From (\dag)  we obtain the exact sequence $0 \to T^{\partial} \to {  G}^{\partial} \to (A^{\flat})^{\partial} \to 0$  of differential algebraic groups, from which, computing points in $K^{diff}$, we finally derive the exact sequence of groups
\newline
(\dag \dag) $0 \to T(C_{K}) \to {  G}^{\partial}(K^{diff}) \to A^{\flat}(C_{K}) \to 0$.

\medskip
\noindent
Let $K_{0}$ be a finitely generated differential field contained in $K$ over which ${ G}$ is defined and such that  $K = (C_{K}.K_{0})^{alg}$. 
We now use the language of ``generic points" over $K_0$, as discussed at the beginning of Section 5.1.

We will first show that we can find a generic point of ${ G}^{\partial}$ over $K_{0}$ which is also a $K^{diff}$-rational point. Let $m = dim( G)$. 
We will use the exact sequence (\dag \dag) above, together with the assumption that $C_{K}$ has 
infinite transcendence degree and $K_{0}$ is finitely generated.  Namely we first choose a point of $A^{\flat}$ which is generic over $K_{0}$, and is $C_{K}$-rational. Let $b \in {G}^{\partial}(K^{diff})$ project to $a$. Now let $c$ be a generic point of $T$ over $K_{0},a,b$, such that $c$ is $C_{K}$-rational. Finally let $d$ be the sum 
 of $b$ and $c$. Then $d\in { G}^{\partial}(K^{diff})$ and it is straightforward to verify that 
$tr.deg(K_{0}(d)/K_{0})$ equals $m$. So $d$ is a generic point of ${ G}^{\partial}$ over $K_{0}$ which is also $K^{diff}$-rational. 

As ${ G}$ is $K$-large, $d\in {  G}^{\partial}(K)$, and so $d$ is in the algebraic closure
in the model-theoretic sense of (the finitely generated) $K_{0}$ together with a finite tuple of elements from $C_{K}$. As every element of ${  G}^{\partial}$ is a product of generic elements, it follows that
(in ${\cal U}$)  ${  G}^{\partial}$ is contained (uniformly) in $acl({\cal C},K_{0})$ and so is
definably isomorphic to the group of ${\cal C}$-points of an algebraic group defined over ${\cal C}$. This implies (see Fact 2.6 of \cite{Kowalski-Pillay} for example) that ${  G}$ is isomorphic to an algebraic group over ${\cal C}$ and hence over $C_{K^{diff}} = C_{K}$. The same is then true of $B$.

Recalling the notations at the beginning of  Subsection 4.2 and in Remark 3.10, we obtain:
\begin{Corollary} Let $K$ be as in the previous proposition. Let $G$ be an almost
semiabelian $D$-group over $K$. Suppose that it is $K$-large. Then $G^{sa}_{0}$ descends to $C_{K}$.
In other words, $({\bf HG})_{0}$ holds of $G$.
\end{Corollary}
\pf  Note that  
${\overline {G^{sa}_{0}}}$ is a $D$-subgroup of a $D$-quotient of $G$. Hence it is also $K$
-large. Now use Proposition 6.1.

\subsection{Differential Galois theory}

\medskip
The point of $K$-largeness is that it allows  
 a Galois theory for equations of the form 
$$\partial\ell n_{G}(-) = a , \quad {\rm where}~  a\in LG(K). \qquad (**)$$ 
Namely, suppose  
$(G,s)$ is a $K$-large  algebraic $D$-group  and 
$\alpha\in G(K^{diff})$ is a solution of (**). Let $F = K(\alpha)$, a differential subfield of $K^{diff}$,  and let $Aut_{\partial}(F/K)$ be the group of automorphisms of the differential field $F$ which fix $K$ pointwise. For $\sigma\in Aut_{\partial}(F/K)$, $\sigma(\alpha)$ is also a solution of (**), so 
$\sigma(\alpha)- \alpha \in Ker(\partial\ell n_{G}) = G^{\partial}$ whereby 
$\sigma(\alpha) = \alpha + \rho_{\sigma}$ for a unique $\rho_{\sigma}\in G^{\partial}(K^{diff})$ and by 
the $K$-largeness assumption in fact $\rho_{\sigma} \in G^{\partial}(K)$. As pointed out in \cite{Pillay}, the map taking $\sigma$ to $\rho_{\sigma}$ establishes an isomorphism between $Aut_{\partial}(F/K)$ and a differential algebraic subgroup of $G^{\partial}(K^{diff})$ which by Fact 2.4 (ii) is of the form $H^{\partial}(K^{diff})$ for $H$ a $D$-subgroup of $G$, defined over $K$. 
Moreover there is a Galois correspondence between differential fields in between $K$ and $F$ and $D$-subgroups of $H$ defined over $K$ (or equivalently, 
by $K$-largeness, over $K^{diff}$).

\medskip

With notation as above, here are some additional remarks  
taken from \cite{Pillay}, to be used below. Working  in ${\cal U}$ and noting that $H^{\partial}$ acts on $G$, we see that the orbit of $\alpha$ under $H^{\partial}$ coincides with its orbit under $Aut_{\partial}({\cal U}/
{\cal C}.K)$ and is a differential algebraic $PHS$ for $H^{\partial}$, defined over $K$.
In particular,  $tr.deg(K(\alpha)/K) = dim(H)$.

\bigskip
We can now 
give the promised Galois theoretic proof of Theorem 1.3 in the $K$-large case.
\begin{Theorem} Let $K$ be algebraically closed and of transcendence degree $1$ over its field of constants $\C$.
Let $G$ be an almost semiabelian $D$-group which is $K$-large. Let $x\in LG(K)$ and $y\in G({\cal U})$ be such that $\partial\ell n_{G}(y) = \partial_{LG}(x)$. Assume $x$ satisfies $({\bf HX})_{K}$. Then $tr.deg(K(y)/K)) = dim(G)$. 
\end{Theorem}
\pf  By Corollary 6.2, $G$ satisfies $({\bf HG})_{0}$, so the theorem follows from Theorem 1.3, 
but the present proof will avoid the difficult Proposition 4.1. 
Suppose the conclusion fails, so it fails for some $y$ in $G(K^{diff})$,  for which we assume $tr.deg(K(y)/K) < dim(G)$. Let the $D$-subgroup $H$ of $G$ be the differential Galois group of $K(y)/K$. As recalled above, $Aut_{\partial}(K(y)/K)$ is isomorphic to $H^{\partial}(K^{diff})$, and since $dim(H) = tr.deg(K(y)/K)$, $H$ is a proper  
 $D$-subgroup of $G$.  
 Let $G' = G/H$. The images $y', x'$ of $y,x$ under the projections $G\to G'$ and $LG \to LG'$ satisfy $\partial\ell n_{G'}(y') = \partial_{LG'}(x')$. Also, by Corollary 6.2 and Lemma 4.2, $({\bf HX})_{K}$ is valid for $x'$ and $LG'$. But now, as the orbit of $y$ under $H^{\partial}$ was defined over $K$, $y'$ is $K$-rational, i.e. in $G'(K)$. Furthermore, for any $D$-quotient $G''$ of $G'$ defined over $K$, these points $y',x'$ project to rational points $y'', x''$ in $G''(K), LG''(K)$ still satisfying $\partial\ell n_{G''}(y'') = \partial_{LG''}(x'')$ and $({\bf HX})_{K}$. Hence we may assume that $G'$ has no proper connected $D$-subgroups. In particular the maximal connected isoconstant  
  $D$-subgroup of $G'$ is either $G'$ itself or $0$. Proposition 4.3 or 4.4 will then give a contradiction.

\vfil \eject

\centerline{\bf APPENDIX}

\bigskip
\centerline{\bf EXPONENTIALS ON ALGEBRAIC $D$-GROUPS}

\appendix

\bigskip
\noindent
By ``the text", we mean the main body of the present article.

\section{Appendix - Setting}

Let  $S$ be a  smooth algebraic  curve over $\C$. In this  appendix, we denote by $K = \C( S)$ the field of rational functions on $S$, {\it not} its algebraic closure. Sometimes, we may withdraw a finite set of points from $S$, but still denote by $S$ the resulting affine curve.   We write $S^{an}$  for the Riemann surface attached to $S(\C)$. Finally, we fix a nowhere vanishing vector field $\partial \in H^0(S, TS)$ on $S$, which we identify with a derivation of $K$, with constant subfield $\C$.

\medskip

We start with a commutative algebraic group $G/K$,   geometrically connected and with split maximal torus. Shrinking $S$ if necessary, we fix a connected group scheme  $\pi : {\bf G} \rightarrow S$ extending $G$ over $S$, all of whose fibers 
 have the same toric and unipotent ranks.  
 We denote by $\bf e$ its $0$-section, and by $L{\bf G}$ the pull-back ${\bf e}^*(T_{{\bf G}/S})$ of the relative tangent bundle of ${\bf G}$ over $S$. In other words, $\bf G$ is an algebraic family $\{Ê{\bf G}_t, t \in S \}$ of commutative algebraic groups over $\bf C$, parametrized by $S$, and $L{\bf G}$ is the algebraic family of their tangent spaces $L{\bf G}_t$ at the origin. At the generic point of $S$, we have the algebraic group $G/K$ with (relative) tangent bundle $T_{G/K} \simeq G \times  LG$; this is denoted by $T(G)$ 
in \S 2 of the text.

\begin{Remark}: We will also need to consider $\bf G$ as an analytic family  ${\bf G}^{an}$ of complex Lie groups over the Riemann surface $S^{an}$. We will drop the 
 exponents $^{an}$ when the context is clear. Moreover, for notational ease, several results below are written at the generic point of $S$, but actually extend to $S$, i.e. can  be ``bold-faced". We can then ``analyticize" them, i.e. add $^{an}$ both on the base  and on the fiber spaces under consideration.
\end{Remark}

\medskip

The (total) tangent bundle $T{\bf G}$ of $\bf G$ sits in an exact sequence
$$0 \rightarrow T_{{\bf G}/S} \rightarrow T{\bf G} \rightarrow \pi^*(TS) \rightarrow 0$$
of vector bundles over $\bf G$, and is also a group scheme over $TS$. When $t$ runs through $S$,  its fibers $(T{\bf G})_{(t, \partial_t)}$ yield a subgroup scheme $T_\partial {\bf G}$ over $S$, whose generic fiber is called the twisted tangent bundle $T_\partial G/K$. A section $\bf y$ of ${\bf G}/S$ provides a section $d{\bf y}$ of  $T{\bf G}/TS$, hence a section $d{\bf y}(\partial)$ of $T_\partial{\bf G}/S$, which in accordance with the text, we denote by  $(y, \partial y) \in T_\partial G(K)$, or sometimes just $\partial y$, at the generic point of $S$. 
 Viewed over $K$, $T_\partial G$ is a group extension of $G$ by $LG$ (in  particular, there is a canonical  identification of $L G$ with the fiber above $e$ of $T_\partial G$, cf. \S 2.2 of the text, and  \cite{Pillay}, \S 2). The  zero section of $T_\partial {\bf G}$ is $d \bf e(\partial)$, written $(e, \partial e) \in T_\partial G$ at the generic point of $S$.

\medskip
Viewed over $G$, $T_\partial G$  is a torsor under $T_{G/K}$, and as such, is described by a class in $H^1(G, T_{G/K})$. Assume now that $\pi$ is proper. We can then consider the Kodaira--Spencer map $\kappa: T_tS \rightarrow H^1({\bf G}_t, T_{{\bf G}_t/\C(t)})$ attached to the abelian scheme ${\bf G}/S$ at the generic point $t$ of $S$. By definition, its value $\kappa(\partial)$ at $\partial$ is the class of the torsor $T_\partial G$. Actually, properness is not required to carry out this construction. In the proper case, it is classical 
 that  $\kappa(\partial)$ vanishes if and only if  the abelian variety $G/K$ descends to $\C$.  See Lemma 3.4.(iii) of the text for the semiabelian case, 
 and Lemma 3.4.(ii) for counterexamples in the general case.

\section{The functors $T_\partial, L$ on algebraic groups}

As explained in \S 2.3 of the text, $T_\partial$ and $ L$ are functorial:  given a morphism of algebraic groups $f : G_1 \rightarrow G_2$, we have the twisted differential $T_\partial f : T_\partial G_1 \rightarrow T_\partial G_2$, with a commutative diagram over $f : G_1 \rightarrow G_2$, and the (vertical) differential at $e_1$ of $f$, namely $Lf : LG_1 \rightarrow LG_2$.   All this can be boldfaced, i.e. comes from group schemes over $S$, where $L$  
now stands for  $d_{{\bf G}/S, {\bf e}}$.  

\medskip

 We will sometimes apply to the algebraic group $T_{\partial}G/K$ itself (and to its extension over $S$) what we are doing on $G$. The following identifications will be crucial.

 \begin{Lemma} Let ${\bf G}/S$ be a group scheme as above. There is a functorial isomorphism between  the  group schemes $L(T_\partial {\bf G})$ and $T_\partial (L{\bf G})$. More precisely, given $ {\bf  f} : {\bf G}_1 \rightarrow {\bf G}_2$, we can identify $T_\partial(L{\bf   f}) : T_\partial (L{\bf G}_1) \rightarrow T_\partial (L{\bf G}_2)$ with $L(T_\partial {\bf   f}) : L(T_\partial {\bf G}_1) \rightarrow L(T_\partial {\bf G}_2)$.
 
\end{Lemma}

\pf This is a straightforward extension over $S$ of Lemma 2.1 of the text, whose proof can be viewed as the study of  the  functor $T$ itself, i.e. of the total differential $d{\bf   f}$ on $T{\bf G}_1$. In what follows, we often write these identifications only at the generic point of $S$.  
Notice that the formula $LT_\partial = T_\partial L$ is compatible with the identification of $L G$ with  $(T_\partial G)_e$, so that the two $LL$'s which they provide coincide. Also, recall that $LLG$ is  canonically isomorphic to $LG$ (more  
generally, 
 when $V$ is a vectorial group, we always identify $V$ and $LV$, but we sometimes keep to the notation $LV$ to  
remove ambiguities).

\medskip
For later use, we point out  that  
$$T_\partial f(y_1, \partial y_1) = (f(y_1), \partial(f(y_1)))$$
 for any $y_1 \in G_1$, and  that  $T_\partial f$ induces $Lf$ on the fiber $(T_\partial G_1)_{e_1}$, identified with $LG_1$.

\section{Algebraic $D$-groups and $\partial \ell n_{\bf G}$}

As in \S 2.2 of the text, we now assume now  that  the group extension $T_\partial G$ is trivial,  in other words, that its class in $H^1(G, T_{G/K})$ vanishes, and we let $s$ be a homomorphic section of $T_\partial G \rightarrow G$, or equivalently, a vector field on $G$ above $\partial$ such that the corresponding derivation on ${\cal O}_{G}$ respects the group structure of $G$.  We then say that $(G,s)$ is a (commutative)  algebraic $D$-group over $K$, and we denote by ADG the corresponding category. 
The {\it logarithmic derivative}  of $(G,s)$ (which should be indexed by $s$) is then defined by
$$\partial \ell n_G : G \rightarrow (T_\partial G)_e \simeq LG : y \mapsto \partial y - s(y).$$ 
Shrinking $S$ is necessary, we can  extend $s$  to a section ${\bf s}: {\bf G} \rightarrow T_\partial{\bf G}$ over $S$, and we then set:
$$\partial \ell n_{\bf G} : {\bf G} \rightarrow (T_\partial {\bf  G})_e \simeq L{\bf G} : {\bf y} \mapsto d{\bf y}( \partial) - {\bf s}({\bf y}).$$
When $V$ is a vectorial group over $K$ (i.e. ${\bf V}/S$ is a vector bundle), $\partial \ell n_V :  V \rightarrow V$ and $\partial \ell n_{\bf V} :  {\bf V} \rightarrow {\bf V}$ are the contractions with $\partial$ of connections in the usual sense.

\medskip
\noindent
We are going to associate to  the logarithmic derivative $\partial \ell n_G$ two (contracted with $\partial$) connections on the vectorial group $LG$:

$\bullet$ an algebraic one: $\partial_{LG}$ (which is the one of the text itself);

$\bullet$ an analytic one: $exp_{\bf G}^*(\partial \ell n_{\bf G}) = \partial \ell n_{\bf G} \circ exp_{\bf G} $;

\noindent this second one is actually defined on $L{\bf G}^{an}$,  but see Remark G.6 for a formal approach. And we will compare  them, {\it in the case $G$ has a unique structure of ADG},  to:
 
 $\bullet$ the  Gauss-Manin connection  $\nabla_{LG,\partial}$ (again, algebraic).
 
 \medskip
More precisely, withdrawing some points of $S$ if necessary, the first one extends over $S$,  i.e. is the value at the generic point of $S$ of a connection $\partial_{L\bf G}$ on the vector bundle  $ L{\bf G}$, and we can look at it analytically, as one on $L{\bf G}^{an}$. Ditto for the third one,  
one of whose characterizations (see Sections H and I) comes from $L{\bf G}^{an}$. In this analytic context, we will prove on the one hand that
 $$\partial_{L\bf G}   = exp_{\bf G}^*(\partial \ell n_{\bf G}),$$
 and on the other hand that
$$\nabla_{L{\bf G},\partial} = exp_{\bf G}^*(\partial \ell n_{\bf G}),$$
so that $\partial_{L\bf G} = \nabla_{L{\bf G},\partial}$. All this is on $(L{\bf G}^{an})/S^{an}$, but these  are equalities, not just isomorphisms. So, we will finally deduce that for any almost semi-abelian $D$-group $G$:
$$\partial_{LG} = \nabla_{LG, \partial} ~{\rm ~Êon}~ LG.$$
However,  the relation with $exp_{\bf G}^*(\partial \ell n_{\bf G})$ will also be useful for other parts of our paper.

\bigskip
\section{The connection $\partial_{LG}$ on $LG$}

As explained in \S 2.3 of the text, this is easy to define in view of Lemma B.1, but we repeat the argument in order to specify which  shrinking of $S$ may be necessary.  The morphism of algebraic groups $s : G \rightarrow T_\partial G$ has a vertical differential $Ls : LG \rightarrow LT_\partial G$ at the zero section $e$, which, under the identification  $LT_\partial G = T_\partial LG$, can be viewed as a section of $T_\partial LG \rightarrow LG$. We have then defined 
$$\partial_{LG} := \partial \ell n_{(LG, Ls)}$$ 
as the logarithmic derivative of the ADG structure on $LG$ defined by this section, i.e.
$$ \partial_{LG} : LG \rightarrow LLG = LG : x \mapsto \partial x - Ls(x).$$
This is a connection on the vector group $LG/K$. 

\medskip
Now, restricting $S$ if necessary, we may assume that $s$ extends to an ${\cal O}_S$-section ${\bf  s}$ of  $T_\partial {\bf G}  Ê\to {\bf G}$. Then, $Ls$ extends over the same base to an ${\cal O}_S$-homomorphism $L{\bf s} =  d_{{\bf G}/S, {\bf e}}({\bf s}) :  L{\bf G} \rightarrow  LT_\partial  {\bf G}$. With the identifications of Lemma B.1 in mind,  the formula 
$$L{\bf G}Ê\ni {\bf x} \mapsto \partial_{L\bf G} {\bf x}Ê = d{\bf x}(\partial)  
-   L{\bf s}({\bf x}) \in L{\bf G} $$
 then defines (the contraction with $\partial$ of) a connection  on the vector bundle $L{\bf G}$, which coincides with $\partial_{LG}$  at the generic point of $S$.

\medskip
We point out that the description,  given at the end of \S 2.3, of $\partial_{LG}$ as the differential $L \partial \ell n_G$ of $\partial \ell n_G$ at the identity, in the sense of Kolchin's differential algebraic geometry, could also be carried out in the present setting.  

\begin{Remark} Let $V$ be a vectorial subgroup of $G$. If $V$ is a ADG subgroup of $(G,s)$, i.e. if $s(V) \subset T_\partial V \subset T_\partial G$, i.e. if $s$ induces a section $s_V$ of $T_\partial V \rightarrow V$, the differential $Ls_V: LV \simeq V \rightarrow LT_\partial V \simeq T_\partial V$ of $s_V$ can be identified with $s_V$, and we get: 
$$(\partial_{LG})|_{LLV = LV} = (\partial \ell n_G)|_{LV = V}.$$
We  will show in Corollary G.4 below that this relation still holds true when the vectorial subgroup $V$ is ~ {\rm not}~ an ADG subgroup of $G$. 
\end{Remark}

\section{The exponential map on ${\bf G}^{an}$} 

 For each $t \in S$, we can consider the exponential map $exp_{{\bf G}_t}: L{\bf G}_t(\C) \rightarrow {\bf G}_t(\C)$ of the connected Lie group ${\bf G}_t^{an}$ attached to ${\bf G}_t(\C)$. Its kernel ${\cal P}_t $ is the $\Z$-module of periods of ${\bf G}_t^{an}$.   These patch into  an exact sequence of analytic sheaves of abelian groups over $S^{an}$ :
$$0 \rightarrow {\cal P}Ê\rightarrow L{\bf G}^{an} \rightarrow {\bf G}^{an} \rightarrow 0 ~,  \qquad (\dagger)$$ 
whose third arrow  $exp_{{\bf G}^{an}}$ induces  $exp_{{\bf G}_t}$ above each $t$. Its kernel $\cal P$ will be described in Section H.  Following Remark A.1, we will   drop the exponents $^{an}$ when the exponential morphism is concerned. Typically, $exp_{\bf G}$ can only mean $exp_{{\bf G}^{an}}$, with source the analytic vector bundle $L({\bf G}^{an}) = (L{\bf G})^{an}$ over $S^{an}$.

\medskip
Let us now collect some properties of this $S^{an}$-morphism $exp_{\bf G} : L{\bf G}^{an}Ê\rightarrow {\bf G}^{an}$. Writing the group law on $G$ additively, it is characterized by the joint conditions that:

\medskip

 (i) $\forall U \subset S^{an}, \forall {\bf x}_1, {\bf x}_2 \in L{\bf G}^{an}(U), exp_{\bf G}({\bf x}_1 + {\bf x}_2) =  exp_{\bf G}({\bf x}_1) +  exp_{\bf G}({\bf x}_2)$;
 
 \medskip
 
(ii) $Lexp_{\bf G} = id_{L{\bf G}}$ (under the usual identification $LL{\bf G} = L{\bf G}$).
 
 \medskip
 \noindent
Here,  $Lexp_{\bf G} = d_{{\bf G}/S,\bf e} \; exp_{\bf G}$ means the vertical differential of $exp_{\bf G}$ along the  zero section $\bf e$, which is still meaningful in the analytic setting.

 \medskip
 So, we must  repeat the whole of Sections A and  B in the context of analytic sheaves of abelian groups over $S^{an}$. More precisely, the situation is as follows: we have two (algebraic) group schemes ${\bf G}_1, {\bf G}_2$ and an {\it analytic} morphism $\phi : {\bf G}_1^{an} \rightarrow {\bf G}_2^{an}$ over $S^{an}$.  We can then define $T_\partial \phi: T_\partial {\bf G}_1^{an} \rightarrow T_\partial {\bf G}_2^{an}$, $L\phi : L {\bf G}_1^{an} \rightarrow L{\bf G}_2^{an}$,  we can again identify the $S^{an}$-sheaves $T_\partial L{\bf G}_i^{an}$ and $LT_\partial{\bf G}_i^{an}$, and get $LT_\partial \phi = T_\partial L\phi$. A typical example will be given by
 $${\bf G}_1 = L{\bf G}, {\bf G}_2 = {\bf G}, \phi = exp_{\bf G}.$$
 From (i) and (ii), we immediately deduce the well-known property:
 
 \medskip
 
(iii) $\forall \phi: {\bf G}_1^{an} \rightarrow {\bf G}_2^{an}$, we have: $\phi \circ exp_{{\bf G}_1} = exp_{{\bf G}_2} \circ L\phi,$
 
 \medskip
 \noindent
and the important fact that

\medskip

(iv) if ${\bf V}$ is a vectorial group scheme ($=$ vector bundle) over $S$, then $exp_{\bf V} = id_{L{\bf V}}$ in the usual identification $L{\bf V} = {\bf V}$.

\medskip
If $G_1$ is a subgroup of $G_2$, Property (iii) shows that $exp_{{\bf G}_1}$ is the restriction to $L{\bf G}_1$ of $exp_{{\bf G}_2}$. Consider in particular the exponential map $exp_{T_\partial{\bf G}} : LT_\partial {\bf G}Ê\rightarrow T_\partial{\bf G}$ of the group scheme $T_\partial{\bf G}$. Its restriction to the Lie algebra $LL{\bf G} = L{\bf G}$ of the vectorial subgroup ${\bf G}_1 = L{\bf G}$ of ${\bf G}_2 = T_\partial{\bf G}$ is the exponential map $exp_{L{\bf G}}$ of $L{\bf G}$. By Property (iv), we therefore have the (trivial, but crucial !) :

\begin{Lemma}  For any algebraic group $G$, 

{\rm (v)}~ $(exp_{T_\partial{\bf G}})|_{LL{\bf G} = L{\bf G} \subset LT_\partial {\bf G}} = exp_{L{\bf G}} = id_{L{\bf G}}$.
\end{Lemma}
This property should not be confused with (ii). (Recall that we are indexing the exponential maps by the groups, not by their Lie algebras.)

\section{The connection $exp_{\bf G}^*(\partial \ell n_{\bf G})$ on $L{\bf G}^{an}$}

Let  
us
make a preliminary comment on pull-backs of connections in the classical case. Let $V_1, V_2$ be two vector spaces over $K$, let $\nabla_2 : V_2 \rightarrow V_2 \otimes \Omega^1_{K/\C}$ be a connection on $V_2$, and let $f: V_1 \rightarrow V_2$ be a $K$-linear map. In general we cannot define the pull-back $\nabla_1 := f^*(\nabla_2)$ of $\nabla_2$ under $f$, but we {\it can} if $f$ is an isomorphism. Indeed, $\nabla_1 := (f\otimes 1)^{-1} \circ  \nabla_2 \circ  f$ is a connection on $V_1$. Notice that it is the unique connection such that $f: (V_1, \nabla_1) \rightarrow (V_2, \nabla_2)$ is a horizontal morphism.

\medskip
Let now $f: G_1\rightarrow G_2$ be a morphism of commutative algebraic groups over $K$, and let $s_2 : G_2 \rightarrow T_\partial G_2$ be an ADG structure on $G_2$. In general, we cannot define the pull-back of $s_2$ nor of $\partial \ell n_{G_2}$ under $f$, but we can if $f$ is an {\it isogeny} (i.e. a finite covering).  Indeed, for any $y_1 \in G_1$ with $y_2 = f(y_1)$, $T_\partial f$ then induces an isomorphism on the fibers $(T_\partial G_1)_{y_1} \rightarrow (T_\partial G_2)_{y_2}$, and we may set $s_1(y_1) = \big( (T_\partial f)_{y_1}\big)^{-1} (s_2(y_2)): G_1 \rightarrow T_\partial G_1$, thereby defining the unique logarithmic derivative $\partial \ell n_{G_1} := f^*(\partial \ell n_{G_2}): G_1 \rightarrow LG_1$ on $G_1$ such that $f$ is horizontal, in the sense that 
$$Lf \circ \partial \ell n_{G_1} = \partial \ell n_{G_2} \circ f.$$
Indeed, $(T_\partial f)(y_1, \partial y_1) = (f(y_1), \partial(f(y_1))$, and  $L$ is the restriction of $T_\partial$ above the zero section, so that $(Lf)(\partial y_1 - s_1(y_1)) = (T_\partial f) ((y_1, \partial y_1) - (y_1, s_1(y_1)) = (T_\partial f) (y_1, \partial y_1)- (T_\partial f)(y_1, s_1(y_1)) = (y_2, \partial y_2) - (y_2 , s_2(y_2)) = \partial (y_2) - s_2(y_2) \in LG_2$.

\medskip
Let us note for the record that the last paragraph also shows that if $G_1$, $G_2$ are algebraic $D$-groups and $f:G_1 \to G_2$ is a $D$-homomorphism (homomorphism of algebraic $D$-groups) then 
$$Lf \circ \partial \ell n_{G_1} = \partial \ell n_{G_2} \circ f.$$

\medskip
The construction above extends word for word over  $S$, and can therefore be boldfaced and analyticized.  Now, more generally, we can consider an analytic morphim $\phi:  {\bf G}_1^{an} \rightarrow {\bf G}_2^{an}$. As soon as $\phi$ is a covering (not necessarily finite, but with discrete kernel), the same construction applies, and 
in parallel with $s_1$,  we obtain an analytic section $\sigma_1 : {\bf G}_1^{an}Ê\rightarrow T_\partial {\bf G}_1^{an}$, hence an analytic logarithmic derivative 
$$\partial \ell n_{{\bf G}_1^{an}} := \phi^*(\partial \ell n_{{\bf G}_2}): {\bf G}_1^{an} \rightarrow L{\bf G}_1^{an}$$
 on ${\bf G}_1^{an}$ such that 
$L \phi \circ \partial \ell n_{{\bf G}_1^{an}} = \partial \ell n_{{\bf G}_2^{an}} \circ \phi.$

\medskip
Since $exp_{\bf G}$ is  
an analytic
 covering, we may apply the latter construction to the situation ${\bf G}_1 = L{\bf G}, {\bf G}_2 = {\bf G}, \phi = exp_{\bf G}$. We thereby  obtain (the contraction with $\partial$ of) an analytic connection on $L{\bf G}^{an}$  :
$$exp_{\bf G}^*(\partial \ell n_{\bf G}) : L{\bf G}^{an} \rightarrow  L{\bf G}^{an},$$
under the usual identification of the vector bundles $LL{\bf G}  = L{\bf G}$ over $S$. 

\medskip

The horizontality of $\phi$ (worked out  in the algebraic setting of  $f$ in the previous section), applied to the present  case $\phi = exp_{\bf G}$, gives 
$$ L exp_G  \circ exp_{\bf G}^*(\partial \ell n_{\bf G}) = \partial \ell n_{\bf G} \circ exp_{\bf G},$$
and
since $ L exp_{\bf G} = id_{L \bf G}$ by Property (ii), we eventually have

\begin{Lemma} 
$exp_{\bf G}^*(\partial \ell n_{\bf G})  = \partial \ell n_{\bf G} \circ exp_{\bf G}.$
\end{Lemma}

\section{$ exp_{\bf G}^*(\partial \ell n_{\bf G}) =  \partial_{L \bf G} $}

Let $G$ be an arbitrary commutative algebraic $D$-group.  So, a section $s: G \rightarrow T_\partial G$ is fixed, and  extended  to $s: {\bf G} \rightarrow T_\partial{\bf G}$ as usual.  Given a local section ${\bf x} \in L{\bf G}^{an}(U)$,  we proceed to compute $exp_{\bf G}^*(\partial \ell n_{\bf G})({\bf x})  = \partial \ell n_{\bf G} \circ exp_{\bf G}({\bf x})$, and will show:

\begin{Proposition} Let $G/K$ be an algebraic $D$-group. Then, $exp_{\bf G}^*(\partial \ell n_{\bf G}) = \partial_{L \bf G}.$
\end{Proposition}
 However,  we will not  bold-face everything in this section (not speaking of the already dropped exponents $^{an}$), in the  hope that the context is clear. 

\medskip
It is at this point that we need the exponential map of the algebraic group $T_\partial G$ itself  (more properly, the exponential morphism of  the analytic family $T_\partial {\bf G}^{an}$). So, we have $exp_{T_\partial \bf G}: LT_\partial {\bf G}^{an} \rightarrow T_\partial {\bf G}^{an}$, abbreviated as $exp_{T_\partial G}$. By the analytic version of  Section B, we may identify its source with $T_\partial L{\bf G}^{an}$.  
Our proof of Proposition G.1 will in particular make use of Lemmas 
E.1 and
F.1.

\medskip
By definition, 
$$\partial_{LG}(x) = \partial x - Ls(x),$$
while
$$\partial \ell n_G (exp_G(x)) = \partial(exp_G(x)) - s(exp_G(x)),$$
By Property (iii), $s(exp_G(x)) = exp_{T_\partial G}(Ls(x)) \in T_\partial G$, and we claim (see Lemma G.3 below) that $\partial(exp_G(x)) = exp_{T_\partial G}(\partial x) \in T_\partial G$. So,
$$\partial \ell n_G (exp_G(x)) = exp_{T_\partial G}(\partial x)  - exp_{T_\partial G}(Ls(x)) =   exp_{T_\partial G}\big(\partial x - Ls(x)\big).$$
But $\partial x - Ls(x) = (x, \partial x) - (x, Ls(x))$ is a point of $ T_\partial L G = LT_\partial G $ lying in the Lie algebra  $LLG$ of the vectorial subgroup $LG$ of $T_\partial G$. By  Property (v) in Lemma E.1, we therefore have: 
$$ exp_{T_\partial G}\big(\partial x - Ls(x)\big) = id_{LG}\big(\partial x - Ls(x)\big) =  \partial x - Ls(x) = \partial_{LG}(x) ,$$
and the proposition is established.

\begin{Remark} This computation, where we lift questions on $G$ to $T_\partial G$, and  then use the vectorial properties of its subgroup $LG$, is reminiscent of methods from universal vectorial extensions. But notice that our algebraic $D$-group $G$ is already (essentially) a universal extension. 
\end{Remark}
 
 \medskip
It remains to show that
\begin{Lemma} For any algebraic group $G$, and any $x$ in $LG$, we have : $\partial(exp_G(x)) = exp_{T_\partial G}(\partial x) \in T_\partial G$.
\end{Lemma}

\pf We already know (analytic version of the last relation in Section B) that $\partial(exp_G(x))  = (T_\partial exp_G)(\partial x)$. We must now prove that 
$$T_\partial exp_G = exp_{T_\partial G}.$$
Since $T_\partial exp_G$ is a group morphism, we are reduced, in view of the properties  (i), (ii)   
characterizing the exponential morphism $exp_{T_\partial G}$ of $T_\partial G$, to check that
$$LT_\partial exp_G = id_{LT_\partial G}.$$
By the analytic version of Lemma B.1, $LT_\partial exp_G  = T_\partial Lexp_G$, and $Lexp_G = id_{LG}$ by Property (ii) for $exp_G$. So, the LHS is $T_\partial id_{LG}$. Now, the RHS is $id_{LT_\partial G} =  id_{T_\partial LG} $, and the requested identity
$$ T_\partial id_{LG} = id_{T_\partial LG}$$
actually holds true for any algebraic group, not necessarily of the type $LG$.

\begin{Corollary} Let $V/K$ be a arbitrary vectorial subgroup of $G$, identified with its Lie algebra $LV \subset LG$, and with $LLV \subset LLG$. Then
$$(\partial_{LG})|_V = (\partial \ell n_G)|_V.$$
In particular, $V$ is an ADG subgroup of $G$ if and only if it is a $\partial_{LG}$-submodule of $LG$.
\end{Corollary}
\pf In view of Property (iv) of the exponential map $exp_V = exp_G|_{LV}$,  we deduce from  Proposition G.1 that indeed, 
$$\partial_{LG} (v) = \partial \ell n_G(exp_G(v)) = \partial \ell n_G(v).$$ 
for any $v \in V$. Notice that contrary to Remark D.1, this formula was not immediately clear. In our paper, it will replace the role of Remark 1.4 of \cite{Buium}, p. 64. It  also clarifies the role of Hypothesis ($\bf H$) at the end of Section I of this Appendix.  

\medskip

The last sentence of Corollary G.4 can be stated in greater generality, as follows. 
\begin{Corollary} Let  $H/K$ be a connected algebraic subgroup of $G$, with Lie algebra $LH$. Then, $H$ is an ADG subgroup of $G$ if and only if $LH$ is a $\partial_{LG}$-submodule of $LG$.
\end{Corollary}
\pf Left to right, which is a special case of Lemma 2.5 of the text,  is clear from  the definition of $\partial_{LG}$ with sections. For right to left, notice that $exp_H = exp_G|_{LH}$, so that for any local section ${\bf y} := exp_{\bf G}({\bf x}), {\bf x} \in {L\bf H}$, of $\bf H$ over a small disk $U$ in $S$, $\partial \ell n_{\bf G}({\bf y}) = \partial_{L\bf G}({\bf x}) \in  L{\bf H}(U)$. Therefore, for any $y \in H$, $\partial \ell n_G(y)$ lies in $LH \simeq  (T_\partial H)_e$, and  $s(y) = \partial y - \partial \ell n_G(y) $ does lie in $T_\partial H$.

\medskip
\begin{Remark} In both Sections D and E, we used only the properties of the exponential morphism of the {\rm formal} group scheme $\hat{\bf G} = $ formal completion of $\bf G$ along the zero section $\bf e$. This formal exponential is again entirely characterized by Properties (i) and (ii) of Section D. Consequently, the whole development of these Sections goes through, with ${\bf G}^{an}$ replaced by $\hat{\bf G}$. The complex-analytical properties of $exp_{\bf G}$, reflected in the study of the full exact sequence $(\dagger)$, are needed only in the next sections.
\end{Remark}

\section{Gauss-Manin: $\nabla_{L{\bf G},\partial} = exp_{\bf G}^*(\partial \ell n_{\bf G})$}

Given a connected commutative algebraic group $G/K$, we extended it to a group scheme ${\bf G}/S$, and considered in Section E  the exponential sequence 
$$0 \rightarrow {\cal P}Ê\rightarrow L{\bf G}^{an} \rightarrow {\bf G}^{an} \rightarrow 0 ~  \qquad (\dagger).$$
Recall  now from Section A that we are assuming that    the fibers ${\bf G}_t, t\in S(\C)$, have constant toric and unipotent ranks,  
which amounts to their having {\it constant topological type}. This can always be achieved by removing a finite set of points from $S$.  Under this assumption, the kernel ${\cal P}$ of $(\dagger)$ is a local system over $S^{an}$, dual to the local system  $R^1\pi_*({\Z})$  formed by the Betti cohomology groups of the fibers ${\bf G}_t$. For any open $U \subset S^{an}$, the sections in ${\cal P}(U) \subset L{\bf G}^{an}(U)$ are killed by the exponential morphism $exp_{\bf G}$. In particular, as soon as $G$ is an algebraic $D$-group,   
Lemma F.1 shows that they are horizontal sections of the analytic connection $exp_{\bf G}^*(\partial \ell n_{\bf G}) = \partial \ell n_{\bf G} \circ exp_{\bf G}$.  
In view of Proposition G.1,   the local sections of  $\cal P$ are therefore horizontal  for the connection $\partial_{L{\bf G}}$:
$$\forall U \subset S^{an}, \forall \lambda_U \in {\cal P}(U), ~\partial_{L\bf G}(\lambda_U) = 0.$$

\medskip
Now, assume that  $L{\bf G}^{an}$ is locally generated over ${\cal O}_{S^{an}}$ by $\cal P$. There then exists {\it at most one} connection on $L{\bf G}^{an}$ killing $\cal P$. So, any {\it  connection} on $L{\bf G}^{an}$ killing $\cal P$ {\it will}  coincide with  $exp_{\bf G}^*(\partial \ell n_{\bf G}) = \partial_{L{\bf G}}$. This is the principle, borrowed from \cite{Malgrange}, on which the whole present Section is based.

\medskip
We are now going to define the Gauss-Manin connection $\nabla_{LG}$ on $LG$ when $G$ is an almost semi-abelian $D$-group. The definition is that of  \cite{Brylinski}  when  $G$ is the universal vectorial  extension of a semi-abelian variety, and the general case follows by taking a quotient.  As in the text, 
we will use the following notations. We write $B$ for the maximal semi-abelian quotient of $G$,   and $\tilde B$ for the universal vectorial extension of $B$; so, $B$ is an extension of an abelian variety $A/K$ by a  torus $T$, which we have  assumed to be  split over $K$.  
 Recall from Lemma 3.4.(ii) that $\tilde B$ carries a unique structure of algebraic $D$-group. By definition of an almost semi-abelian  $D$-group (cf. \S 3.1 of the text),
 there  exists a canonical vectorial subgroup $V$ of $\tilde B$, which is an algebraic $D$-subgroup of $\tilde B$, and such that $G = \tilde B /V$. 
We endow  $G$ with the quotient ADG structure, which is actually the unique ADG structure one can put on $G$, see again Lemma 3.4. So, $\partial \ell n_{\tilde B}$ and $\partial \ell n_G$ are well defined. 

\medskip

In these conditions, the topological hypothesis made on the fibers of the group scheme ${\bf G}/S$  implies that $B/K$ can be continued to  an extension  ${ {\bf B}}/S$ of an abelian scheme ${\bf A}/S$ by the constant torus ${\bf T} =  T \times S$. In other words, the one-motive $M = [0 \rightarrow { \bf B}]$ is smooth over $S$ in the sense of Deligne \cite{Deligne}, III, 10.1.10.  
Then,  $[0 \rightarrow \tilde {\bf B}]$ is  the universal vectorial extension of  $M$,  and the vector bundle $L{\tilde {\bf B}}$ is its de Rham realization $T_{dR}(M)$, see  \cite{Deligne}, {\it loc. cit.}, also \cite{Bertapelle}.  Finally, we denote by $\tilde {\cal P} \subset L{\tilde {\bf B}}$ and $\overline {\cal P} \subset L{{\bf B}}$ the kernels of the exponential exact sequences of the sheaves ${\tilde {\bf B}}^{an}$ and ${ {\bf B}}^{an}$ over $S^{an}$. In particular, $\overline {\cal P}$ is the Betti realization $T_{\Z}(M)$ of the one-motive $M$.

\medskip
Since the exponential maps have no kernel on vectorial groups, the local systems ${\cal P}, \overline {\cal P}, \tilde {\cal P}$ are  isomorphic as abstract  $\Z_{S^{an}}$-sheaves. Tensored with ${\cal O}_{S^{an}}$, they all define the same  vector bundle, say $ {\cal P} \otimes {\cal O} _{S^{an}}$, which, as said above,  carries a unique connection $\nabla = id_{\cal P}\otimes  d_{S^{an}/\C}$ (equivalently $\overline \nabla, \tilde \nabla$) relatively to which they are horizontal. Now, by \cite{Brylinski}, Facts 2.2.2.1 and 2.2.2.2, we have:

\begin{Fact}  The natural map
$$ {\tilde {\cal P}}Ê\otimes {\cal O}_{S^{an}}  \rightarrow L{\tilde {\bf B}}^{an}$$
is an isomorphism of vector bundles over $S^{an}$, and there exists
an algebraic connection $\nabla_{L{\tilde {\bf B}}}$ on the ${\cal O}_S$-module $L{\tilde {\bf B}} = T_{dR}(M)$   such that $\tilde \nabla$ and $\nabla_{L{\tilde {\bf B}}^{an}}$ coincide under this isomorphism (i.e. such that $\nabla_{L{\tilde {\bf B}}}$ kills the local sections of $\tilde {\cal P}$). We define  {\rm  the  Gauss-Manin connection $\nabla_{L\tilde B}$ on $L\tilde B$ } as the connection induced by $\nabla_{L{\tilde {\bf B}}}$  at the generic point of $S$.
\end{Fact}
\noindent
This fact reflects the horizontality of the canonical pairing ``integration of forms against cycles" between $T_{\Z}(M) =  \overline {\cal P} \otimes {\cal O} _{S^{an}}$, endowed with $\overline \nabla$,  and the vector bundle $ H^1_{dR}({ {\bf B}}/S) \simeq T_{dR}(M)^{*}$ formed by the de Rham cohomology groups of the fibers of the  semi-abelian scheme $ {\bf B}/S$, endowed with the algebraic Gauss-Manin connection, as originally defined in \cite{Katz-Oda}.  For a slightly different presentation, see  \cite{Andre} and \cite{Steenbrink-Zucker}.

\begin{Fact}Ê When $B = A$ is an abelian variety,  $\nabla_{L\tilde {\bf A}}$ reduces to the dual of the classical Gauss-Manin connection on the de Rham cohomology group $H^{1}_{dR}({\bf A}/S)$, 
 as described in  \cite{Deligne}, II,  \cite{Coleman}. Notice that Part (ii) of the Corollary below, restricted to the study of $\nabla_{L\tilde {\bf A}}$, is the  only  property of Gauss-Manin connections  needed in the text (see \S 4.2, (I), (II)).
\end{Fact}
\noindent
NB:  in this proper case $B = A$,   the connection $\nabla_p$ described by Buium in  \cite{Buium}, Chapter 3.1,  Remark 1.4, p. 64, coincides with the above $\nabla_{{L\tilde A}, \partial}$.  
For an extension to the general case, see   \cite{Buium}, Chap. 3.2, Theorem 2.2.(3).

\begin{Fact}ÊIn the general case,  $\nabla_{L{\tilde {\bf B}}} $ is an extension, in the category of $D$-modules over $S$, of $\nabla_{L\tilde {\bf A}}$ by  $\nabla_{L {\bf T}}$ (the latter one is a direct sum of copies of the trivial $D$-module $({\cal O}_S, d)$). In other words, $L{\bf T}$ is stable under  $\nabla_{L{\tilde {\bf B}}} $. 
\end{Fact}
\noindent
This merely means that the Gauss-Manin connection respects the weight filtration of the smooth one-motive $M$.  This standard fact from \cite{Brylinski} 2.2.2.1, reflected on each fiber $\tilde {\bf B}_t$ by \cite{Deligne}, III.10.1.8, can in fact  also be deduced from Corollary 3.7 of the text, combined with (the easy side of) Corollary G.5 above.

\medskip
By the definition in Fact H.1, $\nabla_{L{\tilde {\bf B}}}$ kills the local system ${\cal \tilde P} \subset L{\tilde {\bf B}^{an}}$. The principle recalled earlier therefore implies:
$$\nabla_{L{\tilde {\bf B}}, \partial} = exp_{\tilde {\bf B}}^* (\partial \ell n_{\tilde {\bf B}}).$$

Finally, the vectorial subgroup $V$ of $\tilde B$ such that $G = \tilde B/V$ is by hypothesis a $D$-subgroup of $G$. By  Remark D.1 (or  the easy side of Corollary G.5), $V = LV \subset L\tilde B$ is stable under the connection $\partial_{L\tilde B}$, i.e. under $exp^*_{\tilde B}(\partial \ell n_{\tilde B})$ (see Section G), i.e. under $\nabla_{L\tilde B}$ (by what has just been proved). We may therefore {\it define} the algebraic Gauss-Manin connection $\nabla_{LG}$ on $LG = L\tilde B/V$ as the quotient connection induced by $\nabla_{L\tilde B}$ (and everything can be bold-faced). Since $\nabla_{L\bf G}$  kills the image $\cal P $ of $\tilde{\cal P}$ in $L{\bf G}^{an}$, and since $\cal P$ still generates $L{\bf G}^{an}$ locally, the principle again gives:
\begin{Proposition} Let $G/K$ be an almost semi-abelian $D$-group. Then, $ exp_{\bf G}^*(\partial \ell n_{\bf G}) = \nabla_{L{\bf G}, \partial}$ on $L{\bf G}^{an}$.
\end{Proposition}
And at long last:
\begin{Corollary} Let $G/K$ be an almost semi-abelian $D$-group,  let $B$ be its maximal semi-abelian quotient, let $T$ be its toric part, and let $A_0$ be the $K/\C$-trace of its maximal abelian quotient $A$. Then, 

i) $ \partial_{LG}$ coincides with the Gauss-Manin connection $\nabla_{LG, \partial}$;

ii) assume that $B = A \times T$. Then, 

\noindent
(I) the connection  $ \partial_{LG}$ on $LG$ is semi-simple;

\noindent
(II)  its $K$-rational horizontal vectors lie in the image $LG_0(\C)$ of $L\tilde A_0(\C) \times LT(\C)$ in $LG$.

iii) in all cases, $T$ is an ADG subgroup of $G$. 
\end{Corollary}

\pf Combining   Propositions  G.1 and H.4,  we get (i). Assertion (ii), which only requires Fact H.2, then follows from the semisimplicity theorem of  \cite{Deligne}, II, Th\'eor\`eme 4.2.6, and from  the theorem of the fixed part of {\it loc. cit.}, Corollaire 4.1.2. 
  As for (iii),  it of course directly follows from Corollary 3.7 of the text, as  pointed out in the proof of Lemma 3.13.(i). But notice that  it also follows from Fact H.3, together with the less obvious side of  Corollary G.5. 

\begin{Remark} In the general case where $B$ is a non-isotrivial extension, the first part of Assertion (ii) of this Corollary does {\rm not} hold, i.e. $\nabla_{LG}$ is not semi-simple : this can be deduced from the version of Manin's theorem of the kernel  studied in Section J, applied to the dual of $A$. The theorem of the fixed part of  \cite{Steenbrink-Zucker} provides an analogue of  the second part of  (ii),  as illustrated by Claim III of \S 5.3 of the text. \end{Remark} 
 
\section{Almost semi-abelian $D$-groups: an analytic approach.}

We here give a different approach to Section H, which still relies on Fact H.1  to deal with universal extensions (and marginally on Fact H.3), but which may be nicer as it defines $\nabla_{LG}$ directly on $LG$, with no prior analysis of $G$. Let $G/K$ be {\it any} commutative algebraic group, not necessarily an ADG. By Chevalley's theorem, it is a vectorial extension of its maximal semi-abelian quotient $B$. So, $G$ is a push-out of the universal vectorial extension $\tilde B$ of $B$, but not necessarily a quotient of $\tilde B$.

\medskip
 Consider the natural map
$$j : {\cal P} \otimes {\cal O}_{S^{an}} \rightarrow L{\bf G}^{an}.$$
We know that  ${\cal P} \otimes  {\cal O}_{S^{an}}$ carries a unique connection $\nabla$ killing the local system $\cal P$. By Fact H.1 of the previous section, combined with  the isomorphism $\cal P \simeq \tilde {\cal P}$, we know that $\nabla \simeq \tilde \nabla$  ``is" the algebraic Gauss-Manin connection $\nabla_{L{\tilde {\bf B}}}$ on $L{\tilde {\bf B}}$ (which we again denote by $\nabla$ below). The discussion now goes as follows.

\medskip
\noindent
$\bullet$ Assume that  {\it $j$ is an isomorphism}, i.e. both  that $\underline{\Z}_S$-linearly independent periods are ${\cal O}_S$-linearly independent  and  that $L{\bf G}^{an}$ is locally generated by $\cal P$. Then, $G = \tilde B$, and we just set $\nabla_{L G} :=  \nabla$. Notice, as in Malgrange's lecture notes  \cite{Malgrange}, that this hypothesis holds if and only if the fibers ${\bf G}_t$ of ${\bf G}/S$ are all analytically isomorphic to $\C^n/\Z^n \simeq (\C)^{* n}$, i.e. that analytically, they are all isomorphic to a torus $(\G_m)^n$. Of course, they will not be so algebraically  if $B$ has a non-trivial abelian part  (as in  Serre's classical counterexample).

\medskip
\noindent
$\bullet$ Now, merely assume that {\it $j$ is surjective}, i.e. that $L{\bf G}^{an}$ is locally generated by $\cal P$, and consider the kernel ${\cal V} = Ker(j)$ of $j$. This ${\cal O}_{S^{an}}$-module is in general not stable under $\nabla$. For instance, if $G = A$ is proper, $ {\cal V}$ is the first step of the Hodge filtration, cf. \cite{Deligne}, II. 4.4.2 (also III. 10.1.3.1), and $\nabla({\cal V}) \subset {\cal V}$ if and only if $A$ is constant. But assume further that $\cal V$ {\it is} a sub-$\nabla$-module, i.e. that

\medskip
\noindent
$\bullet  \qquad \qquad ({\bf H}) \qquad  Im( j) = L{\bf G}^{an} ~{\rm and} ~  \nabla(Ker(j) ) \subset Ker(j).$

\medskip
\noindent
Then, we can define $\nabla_{L{\bf G}^{an}}$ as the quotient connection which $\nabla$ induces on $L{\bf G}^{an} = ({\cal P} \otimes {\cal O}_{S^{an}})/{\cal V} 
 = L{\tilde{\bf B}}^{an}/{\cal V}$. We do not  know yet that $\nabla_{L{\bf G}^{an}}$ is algebraic.

\medskip
Let us now show that it is,  thereby allowing us to write it $\nabla_{LG}$ at the generic point of $S$ (and finally, to call it the Gauss-Manin connection on $LG$). More precisely:
\begin{Proposition} i) under $(\bf H)$, the connection $\nabla_{L{\bf G}^{an}}$ comes from an algebraic connection $\nabla_{L{\bf G}}$ on $L{\bf G}$;

ii) $(\bf H)$ holds only if (and if) $G$ is an almost semi-abelian $D$-group.
\end{Proposition}
\pf  i)  From Fact H.3, and the well-known fuchsianity of the Gauss-Manin connection in the proper case as in Fact H.2 , we deduce that $\nabla$ too is a fuchsian connection. The $\nabla$-stable  ${\cal O}_{S^{an}}$-module $\cal V$ is therefore algebraic, i.e. of the form ${\bf V}^{an}$, for some ${\cal O}_S$-submodule ${\bf V}$ of 
$L\tilde {\bf B}$. So, the quotient connection is indeed algebraic. Let us temporarily denote by $V'$ the $K$-vector subspace that $\bf V$ defines at the generic point of $S$, so that $LG = L\tilde B/V'$.  
\medskip

ii) Firstly, the surjectivity of $j$ implies that $G$ admits no $\G_a$ factor, i.e. is a quotient  $\tilde B/V$ (rather than just a push-out) of $\tilde B$ by some vectorial subgroup $V/K$ of $\tilde B$. From the relation $LG = L{\tilde B}/ V'$, we deduce  that $LV = V'$. In particular, $LV$ is stable under $\nabla = \nabla_{L\tilde B} $, which is equal to $  \partial _{L\tilde B}$ by the main principle of the previous section. Now, the less immediate part of Corollary G.5 (in fact, here, of G.4) implies that $V$ is an algebraic $D$-subgroup of $\tilde B$. Therefore, $G$ is indeed an almost semi-abelian $D$-group.

\medskip
The converse implication is easier, see the previous section.  Our point here is that {\it the algebraic notion of an almost semi-abelian $D$-group is entirely described by the  
 analytic hypothesis  $(\bf H)$,} whose first assumption can be viewed as a transcendental analogue of Lemma 3.1.(iii) of the text.

\bigskip

\bigskip
\centerline{\bf CONCLUSION}

\bigskip
Given an algebraic $D$-group $G/K$, with logarithmic derivative $\partial \ell n_G$, we have  
constructed four connections  (contracted with $\partial$) on its Lie algebra $LG$, and they all coincide:

\bigskip
$\bullet$  the purely algebraic  $\partial_{LG}$, as used in the text;

$\bullet $ the differential-algebraic $L\partial \ell n_{G}$ (not used);

$\bullet$ the analytic  (and  actually formal) $\exp_G^*(\partial \ell n_G) = \partial \ell n_G \circ exp_G$;

$\bullet$ the purely algebraic $\nabla_{LG, \partial}$ (if $G$ is an almost semi-abelian $D$-group).

\medskip
\noindent
In what follows, as well as in the text itself, we identify them, using only the notation $\partial_{LG}$, 
which  we call the logarithmic derivative of $LG$.  

\bigskip
The next Sections of the Appendix are of a different nature. As seen in the text, our main theorem 1.3 relies in an essential way on the Manin-Coleman-Chai  theorem of the kernel.  Although this result is  well documented in the literature (see \cite{Manin}, \cite{Coleman-Manin}, \cite{Chai}), the topic is delicate (see \cite{Coleman-Manin}, \cite{Bertrand-Chai}) and it has seemed useful to provide the reader with a self-contained proof of what we need. We have taken opportunity of this exercise to rewrite the results  in the language of logarithmic derivatives.

 \section{Manin's theorem}

 \bigskip
In this section, we present a weak version of  Manin's theorem of the kernel, in preparation for the next section on Chai's sharpening, as needed in the last step  of both proofs of our main theorem, see Proposition 4.4,  \S 5.1 and \S 6 of the text. The proof given here is essentially the analytic one of Coleman \cite{Coleman-Manin}, Thm. 1.4.3, which, at this point,  is close to  Manin's initial proof \cite{Manin}.
 
  \medskip 
 
 From now on, we restrict to the proper case of an abelian scheme over $S$, as in Fact H.2 of Section H. So, an abelian variety $A/K$ is given and extended to an abelian scheme   $\pi : {\bf A} \rightarrow S$. We have  its $0$-section  $\bf e$, and  the pull-back $L{\bf A} = {\bf e}^*(T_{{\bf A}|S})$ of the relative tangent bundle of ${\bf A}$ over $S$.   As in the introduction to the paper,  we let $(A_0, \tau)$ be the $K/\C$ trace of $A/K$. After base change to a finite cover of $S$,
 we may assume that it also is the $K^{alg}/\C$-trace of $A/K^{alg}$, and that $\tau$ is an embedding. We denote by ${\bf A}_0 := A_0 \times S$ the abelian subscheme of $\bf A$ extending $A_0$.
  
 \medskip
As in Section H, the kernel  of the  exact exponential sequence of analytic sheaves over $S^{an}$ :
$$0 \rightarrow {\cal P} \rightarrow L{\bf A}^{an} \longrightarrow {\bf A}^{an}
\rightarrow 0 , \quad (\dagger_A)$$
is a local system ${\cal P} $ 
  of $\Z$-modules of rank $2g = 2 dimA$ over $S^{an}$, which can be identified with the dual of the local system   $R^1\pi_*({\Z})$ formed by the Betti cohomology groups of the fibers ${\bf A}_t$.

\medskip
Since $\pi$ is proper, any point $y$ of $A(K)$ extends uniquely to a section ${\bf y}Ê\in {\bf A}(S)$, and we will freely use the transition from normal to bold-face characters. On the Lie algebra level, a point $x$ of $LA(K)$ extends to a section of $L\bf A$ only over a Zariski open, but still dense, subset of $S$ (which depends on $x$).

 \medskip
 \noindent
 \begin{Lemma}: let ${\bf x} \in L{\bf A}^{an}(S^{an})$. Assume that ${\bf y} := exp_{\bf A}({\bf x}) \in {\bf A}^{an}(S^{an})$ actually lies in ${\bf A}(S)$. Then, $\bf y$ is infinitely divisible in ${\bf A}(S)$. In particular, there exists a positive integer $d$ such that $d.y$ lies in $A_0(\C)$.
 \end{Lemma}
 
 \pf 
 (Manin-Shafarevich).- The last statement follows from the functional Mordell-Weil theorem, according to which the group $A(K)/A_0(\C)$ is finitely generated, or more directly, from the study of the N\'eron-Tate height $\hat h$ on $A$ attached to an ample divisor, namely:  
 since $y$ is infinitely divisible in $A(K)$, its height $\hat h(y)$ vanishes, while the divisible hull $ A_0(\C)^{div}$ of $A_0(\C)$ in $A(K^{alg})$ is precisely the set of points with zero height.  We now prove the first statement.

 \medskip
 
Let $m$ be an arbitrary positive integer. Then, ${1\over m}{\bf x}$ is again a section of $L{\bf A}^{an}$ over $S^{an}$, and we can consider the  analytic section ${\bf y}_m := exp_{\bf A}({1\over m} {\bf x})$ of ${\bf A}^{an}$ over $S^{an}$. Since $m{\bf y}_m = exp_{\bf A} (m. {1\over m} {\bf x}) = {\bf y}$, $\bf y$ is divisible by $m$ in ${\bf A}^{an}(S^{an})$. Moreover, ${\bf y}_m$ has moderate growth at the points at infinity   of $S$, since its coordinates are algebraic over $K = \C(S)$. 
So, ${\bf y}_m$ actually lies in ${\bf A}(S)$, and $\bf y$ is infinitely divisible in ${\bf A}(S)$, as was to be shown.

\bigskip
 Let now $\tilde A$ be the universal vectorial extension of $A$, endowed with its unique ADG structure, and let $\partial \ell n_{\tilde A}, \partial_{L\tilde A}$ be the corresponding logarithmic derivatives on $\tilde A, L\tilde A$. As in the text, we denote by $W_A \simeq H^1(A, {\cal O}_A)^*$ the maximal vectorial subgroup of $\tilde A$; it is defined over $K$, and may be viewed as vector subspace of $L\tilde A$, but is usually not a $D$-submodule of $L\tilde A$.  
 
  \medskip
 \noindent
 \begin{Proposition} (Manin-Coleman): let $\tilde x \in L\tilde A(K)$, and let $\tilde y \in \tilde A(K)$ with projection $y \in A(K)$. Assume that   $ \partial \ell n_{\tilde A} \tilde y = \partial_{L\tilde A} \tilde x$. Then,  there exists a positive integer $d$ such that $d. y \in A_0(\C)$; in particular, $\tilde x \in L\tilde A_0(\C) + W_A(K)$. 
\end{Proposition}
\pf  Shrinking $S$ if necessary, we may assume that $\tilde x, \tilde y$ extend to sections 
$\tilde{\bf x}, \tilde{\bf y}$ of $L\tilde{\bf A}, \tilde{\bf A}$. Then,   $ \partial \ell n_{\tilde {\bf A}} \tilde {\bf y} = \partial_{L\tilde {\bf A}} \tilde {\bf x}$.  By the main result of Section G,  $\partial \ell n_{\tilde {\bf A}} (exp_{\tilde {\bf A}}(\tilde {\bf x})) = \partial_{L\tilde{\bf A}} \tilde {\bf x}$, so that 
 $$\partial \ell n_{\tilde {\bf A}} \big(\tilde {\bf y} - exp_{\tilde {\bf A}}(\tilde {\bf x})\big) = 0.$$
 In other words, $\tilde {\bf y} - exp_{\tilde {\bf A}}(\tilde {\bf x})$ lies in the kernel $\tilde A^{\partial}$ of $\partial \ell n_{\tilde A}$, and we could continue by projecting it to a point $y_0$ in the  
 differential algebraic subgroup $A^\sharp$ of $A$ defined in Proposition 3.9 of the text.  But $y_0$ is a priori {\it not} defined over the differential closure $K^{diff}$ of $K$, only over the field of meromorphic functions on $S^{an}$, and we cannot appeal to the  
description of $A^\sharp(K^{diff})$   given at the beginning of \S 6.  So, we go back to complex analysis.
 
 \medskip
 
 By the surjectivity of the sheaf morphism $exp_{\tilde {\bf A}}$ from the analogue 
 $(\dagger)$ 
 of $(\dagger_A)$ for $\tilde {\bf A}$, there exists a covering of $S^{an}$ by disks $U$, and sections $\lambda_U \in L{\tilde {\bf  A}}^{an}(U)$ such that $\tilde {\bf y}|_U = exp_{\tilde {\bf A}}(\lambda_U)$.  Two $\lambda_U, \lambda_{U'}$ differ by an element in the kernel $\tilde {\cal P}(
 U \cap U')$
  of $(\dagger)$, allowing to define a cocycle  $\{\phi_{U,U'} := (\lambda_U - \lambda_{U'})\} \in H^1(S^{an}, \tilde {\cal P})$.  In view of Proposition G.1, $\partial \ell n_{\tilde {\bf A}} (exp_{\tilde {\bf A}}(\lambda_U)) = \partial_{L\tilde{\bf A}} \lambda_U$. The relation $ \partial \ell n_{\tilde {\bf A}} \tilde {\bf y} = \partial_{L\tilde {\bf A}} \tilde {\bf x}$ then implies that $\xi_U :=  \lambda_U - \tilde {\bf x}|_U$ lies in $\C$-vector space of horizontal sections of $ \partial_{L\tilde{\bf A}}$ over $U$. But since $ \partial_{L\tilde{\bf A}} =   \nabla_{L\tilde{\bf A}}$ is the Gauss-Manin connection, this space  $({L\tilde{\bf A}})^{\partial}(U)$ coincides with $(\tilde {\cal P} \otimes \underline \C) (U)$, cf. Facts H.1 and  H.2. So, for each $U$,  
 $$\exists \lambda_U \in L{\tilde {\bf  A}}(U),  \exists \xi_U \in (\tilde {\cal P} \otimes \underline \C) (U) ~ {\rm{such~ that}} ~ \xi_U :=  \lambda_U - \tilde {\bf x}|_U.$$
 Since $\tilde {\bf x}$ is a global section over $S^{an}$ (in fact,  even over $S$),  $\lambda_U - \lambda_{U'} =  \xi_{U} - \xi_{U'}$, and  the cocycle $\phi$ has trivial image under the natural map 
  $$H^1(S^{an}, \tilde {\cal P}) \rightarrow H^1(S^{an}, \tilde {\cal P} \otimes \underline \C) = H^1(S^{an}, \tilde {\cal P}) \otimes \C.$$
  This implies that $\phi$ is a torsion point in $H^1(S^{an}, \tilde {\cal P})$. In other words, on refining the covering $\{U\}$ of $S^{an}$ if necessary, there exist a positive integer $d$ and sections $\xi^{1}_U \in \tilde {\cal P} (U)$ such that on all intersections  $U \cap U'$,
  $$d . \lambda_U - d. \lambda_{U'} = \xi^{1}_U - \xi^{1}_{U'}.$$
  In particular, the sections $\tilde {\bf x}^{1}_U := d. \lambda_U - \xi^{1}_U$ glue into a global section $\tilde {\bf x}^{1}$ in $ L\tilde {\bf  A}^{an}(S^{an})$, and since the $\xi^{1}_U$'s lie in the kernel of $exp_{\tilde {\bf A}}$, we finally obtain:
  $$exp_{\tilde {\bf A}} (\tilde {\bf x}^{1}) =exp_{\tilde {\bf A}} (d. \lambda_U) = d. \tilde {\bf y}.$$ 
Let $  {\bf x}^{1} \in L{\bf A}^{an}(S^{an}),  {\bf y} \in {\bf A}(S)$ be the images of $\tilde {\bf x}^{1}, \tilde {\bf y}$ under the natural projections. Since the exponential morphisms commute with these projections, we obtain $exp_{\bf A}({\bf x}^{1}) = d.{\bf y}$, and  these sections exactly satisfy the hypotheses of the previous Lemma J.1. Consequently, a multiple by a positive integer of $d.y$, hence also of $y$, lies in $A_0(\C)$, and the main part of the Proposition is proved.

\bigskip
Recall the notation $W_A$ above. To prove the last sentence, and connect back with the running hypothesis $({\bf HX})_K$ of the paper, we note:
\begin{Lemma} Let $A/K$ be an abelian variety,  let $F$ be a differential extension of $K$ with $F^\partial = \C$, and let  $\tilde x , \tilde y$ be $F$-rational points on $L\tilde A, \tilde A$, projecting onto points $x, y$ in $LA, A$. Consider the following properties:
  
 i) there exists a positive integer $d$ such that $d. y \in A_0(\C)$;
 
 ii) there exists a positive integer such that $d.\tilde y \in \tilde A_0(\C) + W_A(F)$;
 
 iii) $x$ lies in $LA_0(\C)$
 
 iv) $\tilde x$ lies in $L\tilde A_0(\C) + W_A(F)$.

\noindent
Then, $(i) \Leftrightarrow ii)$ and $(iii) \Leftrightarrow (iv)$. Moreover, if $\partial \ell n_{\tilde A} \tilde  y = \partial_{L\tilde A} \tilde x$ {\bf and} if  $F = K$,  they all hold true.
 \end{Lemma} 
 \pf $ii) \Rightarrow i)$ and $iv) \Rightarrow iii)$ are obvious, since $A = \tilde A/W_A$. For $i) \Rightarrow ii)$, notice that any point $y' = d.y$ in $A_0(\C)$ lifts to a point $\tilde y' \in \tilde A_0(\C)$, which will differ from the given $F$-rational lift $d.\tilde y$  by an $F$-rational point in $W_A$. Ditto for $iii) \Rightarrow iv)$. Finally, we have proved above that the two further assumptions imply $i)$, and it remains to show, say, that 
 under these assumptions,  $ii) \Rightarrow iv)$. Indeed, $ii)$,  written as $d.\tilde y = d.\tilde y_0 + w$, together with  the differential relation, the divisibility of vector spaces,  and Corollary G.4 applied to $W_A$, implies that $\partial_{L\tilde A} \tilde x = \partial \ell n_{\tilde A} \tilde  y_0 + \partial \ell n_{\tilde A} w = \partial_{L\tilde A} (w)$, and $\tilde x - w$ is a  horizontal vector of $\partial_{L\tilde A}$, rational over $F = K$,  hence in $L\tilde A_0(\C)$, in view of Corollary H.5.(ii).

 \section{Chai's sharpening}
Let $A/K$ be an abelian variety, with universal extension $\tilde A$, endowed with its canonical ADG structure.  We again  denote by $W_A \simeq H^1(A, {\cal O}_A)^*$ the maximal vectorial subgroup of $\tilde A$, and by $U_A$ its maximal vectorial $D$-subgroup.  
They are defined over $K$, and can be viewed as vector subspaces of $L\tilde A$.  

As in \S 4.3 of the text, define an almost  abelian $D$-group $G$ as an almost semi-abelian $D$-group with no toric part, i.e. such that $B = A$ is an abelian variety. So, $G$ is a quotient of $\tilde A$ by a $D$-vector subgroup $V  \subset U_A \subset W_A$, and $G$ is endowed with its unique ADG  structure. Assuming for simplicity that $A_0 = \{0\}$, and setting $W_G := W_A/V$, we now present Chai's sharpening \cite{Chai} of Manin's theorem,  in terms of logarithmic derivatives.

\begin{Theorem} (Chai) Let $G/K$ be an almost abelian $D$-group such that the abelian variety $ A$ is traceless. Let $x \in LG(K), y \in G(K)$ satisfy $ \partial \ell n_G y = \partial_{LG}   x$. Then,   $ x \in LW_G(K)$. 
\end{Theorem}
\pf  The  proof given here is essentially the dual of the  proof proposed by Chai at the end of his paper \cite{Chai}.
 Instead of using Cech cohomology as for Manin's theorem, we will describe the argument in terms of Galois cocycles. Since  we need only one disk $U$, even  only one point $t_0$ in $S^{an}$, to define them,  we will not use bold-face letters. We let  the fundamental group $\pi_1 = \pi_1(S(\C), t_0)$ act by analytic continuation on the local sections near $t_0$  of any local system over $S^{an}$.  In parallel with the kernel $\tilde {\cal P}$ of $exp_{\tilde {\bf A}}$, we have the local systems $   {\cal P}, \overline {\cal P}$,  respectively defined as the kernels of the morphisms $exp_{\bf G}$ and $exp_{\bf A}$; we recall from Section H that they are all isomorphic.     
We denote by $(L\tilde A)^\partial \simeq \tilde {\cal P}Ê\otimes \C$, resp. $ (LG)^{\partial}$,  the vector spaces of analytic solutions of $\partial_{L\tilde A}$, resp.  $\partial_{LG}$, near $t_0$  (an exponent $^{an}$ will be added when the context is ambiguous.)  As usual, we identify $V$ and the $D$-submodule $LV$ of $L\tilde A$. We  have $V \subset W_A, A = \tilde A/W_A, G = \tilde A/V, LG = L\tilde A /LV$.

\medskip

Let $\tilde y$ be an arbitrary lift of $y$ to a point in $\tilde A(K)$, regular near $t_0$. We will show that $x$ lifts to a point $\tilde x$ in $L\tilde A(K)$ such that
$$ \partial \ell n_{\tilde A} \tilde y = \partial_{L\tilde A} \tilde x.$$
Since $A_0 = {0}$,  Manin's theorem from the previous section then implies that $\tilde x \in W_A(K)$, and the conclusion $x \in W_G(K) (= LW_G(K))$ will follow by projecting modulo $V$. Notice that two lifts differ by an element of $V(K)$, so that this assertion must be independent of the choice of $\tilde y$, in view of  Corollary G.4; in fact, this corollary shows that its truth  depends only on the projection $\overline y$ of $y$ to $A$, as should be.

\medskip
Let $\tilde \lambda$ be a local section of $L\tilde A^{an}$ such that $\tilde y = exp_{\tilde A}( \tilde \lambda)$, and more generally, let ${\cal P}_{\tilde y}$ be the $\Z_{S^{an}}$-local system  formed by all the determinations of the logarithms of all  the multiples of $\tilde y$ in $\tilde A(K)$. Denote by $\Gamma_{\tilde y}$ the $\Q$-Zariski closure of the image of $\pi_1$, acting on ${\cal P}_{\tilde y}\otimes \Q$. Since $exp_{\tilde A}$ is uniform on $S^{an}$  and $\tilde y $ is $K$-rational,   $\gamma \tilde \lambda - \tilde \lambda$ lies in  $\tilde {\cal P}$ for any $\gamma \in \pi_1$. This expression defines a   cocycle 
$$\chi \in H^1(\Gamma_{\tilde y}, \tilde  {\cal P} \otimes \Q):  \Gamma_{\tilde y} \ni \gamma \mapsto  \chi(\gamma) = \gamma \tilde \lambda - \tilde \lambda \in \tilde {\cal P}\otimes \Q.$$
which represents the class of  ${\cal P}_{\tilde y} \otimes \Q$, viewed as an extension of  the local system  $\Q_{S^{an}}$ by $\tilde {\cal P}\otimes \Q$. 

\medskip
Let  further $ N$ be the kernel of the representation $\tilde {\cal P}Ê\otimes \Q  \subset  {\cal P}_{\tilde y}\otimes \Q$ of $\Gamma_{\tilde y}$.  Since $(L\tilde A)^\partial \simeq \tilde {\cal P}Ê\otimes \C$ and since the connection $\partial_{L\tilde A}$ is fuchsian, its differential Galois group  is  the extension $\Gamma_\C$ to $\C$  of the quotient $ \Gamma = \Gamma_{\tilde y}/N$. And it  is  a reductive group,  in view of the semisimplicity given by Corollary H.5.(ii). On the other hand,  $\chi$ induces on the normal subgroup $N$ of $\Gamma_{\tilde y} $  a   $\Gamma$-equivariant injective homomorphism 
 $$\xi  = \chi_{|N} :  N \hookrightarrow \tilde{ \cal P}\otimes \Q  \subset (L\tilde A)^{\partial}.$$
 We now study  the image $\xi(N) \subset \tilde  { \cal P}\otimes \Q$ of $\xi$. We will first show that it vanishes, and then deduce from the reductivity of $\Gamma_\C$ a construction of the desired point $\tilde x \in L\tilde A(K)$. 
 
\medskip
To check that $\xi(N) = \{0\}$, it suffices to show that it is contained in $(LV)^{\partial}Ê  \subset (L\tilde A)^{\partial}$. Indeed, $LV \simeq  V$ is contained $W_A$, which intersects $\tilde{ \cal P}\otimes \Q$ only at $\{0\}$, since in the projection from $L\tilde A$ to $LA = L\tilde A/W_A$, the local system $\tilde {\cal P}$ maps isomorphically onto $\overline {\cal P}$ (see Section I for the relation with the Hodge filtration). 
Now, to prove that $\xi(N) \subset (LV)^{\partial}$, consider the ADG projections $p : \tilde A \rightarrow G$,  $Lp: L\tilde A \rightarrow LG$,  let  $\lambda = Lp(\tilde \lambda) \in  LG^{an}$, and set 
$$\chi_{LG}(\gamma) := Lp(\chi(\gamma)) = \gamma \lambda - \lambda \in  {\cal P}.$$
 Then, $exp_G(\lambda) = p (exp_{\tilde A}(\tilde \lambda)) = y$, and $\partial \ell n_{G}  y = \partial \ell n_{G} (exp_{G} \lambda) = \partial_{LG}  \lambda.$ But by hypothesis, this is equal to $\partial_{LG} x$, with $x \in LG(K)$, so that $\lambda':= \lambda - x \in (LG)^\partial$ is a horizontal section of $\partial_{LG}$. Since $x$ is $K$-rational, $\gamma(\lambda - \lambda') = \lambda - \lambda'$, and we get 
$$\exists \lambda' \in (LG)^\partial, \forall \gamma \in  \Gamma_{\tilde y}, ~ \chi_{LG}(\gamma) = \gamma \lambda' - \lambda'.$$
(In other words, the image $\chi_{LG}$ of $\chi$ in $H^1( \Gamma_{\tilde y}, p(\tilde {\cal P}) = {\cal P})$ vanishes in $H^1( \Gamma_{\tilde y}, (LG)^\partial = j({\cal P}\otimes \C))$, in the notations of Section I.) But $N$ acts trivially on $(L\tilde A)^{\partial}$, hence on its quotient $(LG)^{\partial}$. Restricting $\chi_{LG}$ to $N$, we therefore obtain: $Lp(\xi(\gamma)) =  \gamma \lambda' - \lambda' = 0$ for all $\gamma \in N$. So   $\xi(N)$ does lie in $Ker(Lp) = LV$, and 
 therefore $N \simeq \xi(N) = 0$.

\medskip

Consequently,  $\Gamma_{\tilde y}$  coincides with $\Gamma$.  In particular, the action of the differential Galois group  $\Gamma_\C$ on $(L\tilde A)^\partial$ lifts to an action on the affine space  of solutions of the inhomogeneous equation $\partial_{L\tilde A}(-) = \partial \ell n_{\tilde A}  \tilde y$.  Since $\Gamma_\C$ is  
a reductive group, the corresponding PHS under  $(L\tilde A)^\partial$ is trivial, and the latter equation must  admit  a $K$-rational solution $\tilde x \in L\tilde A(K)$. We now show that $\tilde x$ satisfies the required conditions. On the one hand, $ \partial \ell n_{\tilde A} \tilde y = \partial_{L\tilde A} \tilde x$ by construction. On the other hand, the projection $x'$ of $\tilde x$ to  $LG$ satisfies $\partial_{LG}(x'-x) = \partial \ell n_G y -  \partial \ell n_G y = 0$; since $\partial_{L\tilde A}$ is semi-simple, $x'-x$ lifts to a $K$-rational point of $(L\tilde A)^\partial$, which  vanishes  in view of  the fixed part theorem in Corollary  H.5.(ii) and our hypothesis that  $A$ is traceless. So, $x' = x$, and $\tilde x$ is a lift of $x$, as required.

\medskip
\begin{Remark} The fact that the D-submodule $LV$ of $L\tilde A$ is contained in $LW_A$ has played a crucial role in the proof that $N= 0$. Chai presents in \cite{Chai} a more general version of his theorem, dealing with arbitrary non-zero $D$-submodules of $L\tilde A$. His proof  
of this more general result
actually requires an additional hypothesis, but the statement can be shown to hold in full generality, cf. \cite{Bertrand-Chai}. 
\end{Remark}

\bigskip
We conclude with  an application to sharp points.  Here, $K^{alg}$ denotes the algebraic closure of the field $K = \C(S)$, endowed with the extension of the  derivation $\partial$. Recall from the beginning of \S 6 of the text on $K^{alg}$-largeness 
 that  the sharp points of an abelian variety $A$ always satisfy $A^\sharp(K^{diff}) = A^\sharp(K^{alg})$. We now turn to a description of  the group  $A^\sharp(K^{alg})$ itself. The exponent $^{div}$ still denotes the divisible hull.
\begin{Corollary} Let $A/K^{alg}$ be an abelian variety and let  $A_0$ be its $K^{alg}/\C$-trace. Then, the differential algebraic group $A^\sharp$  satisfies: $A^\sharp(K^{alg}) = A_0(\C)^{div}$. 
\end{Corollary}

\pf  The result is clear  if $A = A_0$. By a standard argument, it remains to study the case when $A$ is  traceless.   
By Proposition 3.9 of the text,  the projection $\pi$ from $G := \tilde A/U_A$ to $A$ induces an isomorphism of differential algebraic groups between the kernel $G^\partial$ of $\partial \ell n_G$ and  $A^\sharp$. If $y \in G(K^{alg})$ satisfies $\partial \ell n_G y = \partial_{LG} x$ with $x = 0$, the previous theorem K.1 (more properly, the lifting property we actually proved at the level of $\tilde A$), combined with the last sentence of Lemma J.3,  implies that the projection  of $y$ to $A$ is a torsion point.

\vfil \eject

\end{document}